\documentclass[preprint,12pt]{elsarticle}

\usepackage{bm}
\usepackage{subfigure}
\usepackage{caption}
\usepackage{amssymb}
\journal{XXXX}
\begin{document}

\begin{frontmatter}



\title{A fast semi-discrete Kansa method to solve the two-dimensional spatiotemporal fractional diffusion equation}

\author{HongGuang Sun \& Xiaoting Liu}
\address{State Key Laboratory of Hydrology-Water Resources and Hydraulic Engineering, College of Mechanics and Materials, Hohai University, Nanjing 210098, China}
\author{Yong Zhang}
\address{1. College of Mechanics and Materials, Hohai University, Nanjing 210098, China\\
2. Department of Geological Sciences, University of Alabama, Tuscaloosa, AL 35487, USA\\
Corresponding author: yzhang264@ua.edu}
\author{Guofei Pang}
\address{Beijing Computational Science Research Center, Beijing 100084, China}
\author{Rhiannon Garrard}
\address{Department of Geological Sciences, University of Alabama, Tuscaloosa, AL 35487, USA}
\begin{abstract}
Fractional-order diffusion equations (FDEs) extend classical diffusion equations by quantifying anomalous diffusion observed frequently in heterogeneous media. Real-world diffusion can be multi-dimensional, requiring efficient numerical solvers that can handle long-term memory embedded in mass transport. To address this challenge, a semi-discrete Kansa method is developed to approximate the two-dimensional spatiotemporal FDE, where the Kansa approach discretizes the FDE first, the Gauss-Jacobi quadrature rule then solves the corresponding matrix, and finally the Mittag-Leffler function provides an analytical solution for the resultant time-fractional ordinary differential equation.  Numerical experiments are then conducted to check how the accuracy and convergence rate of the numerical solution are affected by the distribution mode and number of spatial discretization nodes. Applications further show that the numerical method can efficiently solve two-dimensional spatiotemporal FDE models with either a continuous or discrete mixing measure.  Hence this study provides an efficient and fast computational method for modeling super-diffusive, sub-diffusive, and mixed diffusive processes in large, two-dimensional domains.
\end{abstract}
\begin{keyword}
Anomalous transport \sep Spatiomteporal FDE \sep Semi-discrete Kansa method \sep
Vector fractional
derivative
%\PACS
\end{keyword}
\end{frontmatter}
%\linenumbers*[1]      % YZ add on 5/04/09
\section{Introduction}
Anomalous or Non-Fickian diffusion is a common feature of contaminant transport in natural, potentially multi-dimensional geologic media such as streams, aquifers, soils, and fractured networks \cite{Amir_Neuman2001,Metzler2000,Berkowitz2000,Benson2000}, where the resulting diffusive contaminant plume can cover large areas via super-diffusion.  Contaminant transport modeling has been the prerequisite for many environmental protection applications, but efficiently predicting real-world mass transport remains a historical challenge for the geophysical and hydrological communities \cite{Zhang2009,Schumer2009}.  The physical heterogeneity of natural geological media with intrinsic scale effects has been well recognized as the major mechanism of anomalous diffusion, but it cannot be measured exhaustively at all relevant scales \cite{Dagan1984,Boso2013}.  It is well-known that traditional second-order diffusion models built upon the standard Fick's law cannot efficiently capture non-Fickian features of mass transport processes \cite{Neuman2009, Berkowitz2002,Meerschaert2002,Sun2009,Huang2008}.  These features include apparent early-arrivals and/or heavy late-time tailing in contaminant breakthrough curves, as well as scale-dependent transport coefficients \cite{Berkowitz2006,Seymour2007,Berkowitz2002}.  The failure of second-order diffusion equations in capturing non-Fickian transport has motivated the application of time and space fractional-order diffusion equations (FDEs) to contaminant transport in porous media \cite{Benson2000,Zhang2009}.

Although the FDE models have been successfully used to capture anomalous diffusion in the last several decades \cite{Zaslavsky2002,Klafter2005,Magin2008}, efficient numerical solvers for the FDEs are still required for engineering applications that involve large, geometrically complex spatial domains and long-time range predictions.  The need for these solvers stems from two reasons.  First, the time fractional derivative in the FDE retains a memory of its entire history.  Compared with the integer-order derivative which is a local operator, the computational cost of the temporal FDE dramatically increases with the total modeling time.  Second, although the global correlation defined by the space fractional derivative (especially the fractional Laplacian operator) enables the FDE to characterize the spatial non-local nature of contaminant transport, it causes a computational challenge in approximating spatial FDEs using the grid-based finite difference method (FDM) or finite element method (FEM) \cite{Roop2006,Podlubny2009,Li2012,Jin2015}.  Because the fractional derivative is a non-local operator, the parameter matrix in the FDM (or the stiffness matrix in the FEM) for FDEs is not banded or sparse (as is the standard diffusion equation), but almost dense \cite{Roop2006}.  This issue challenges the accuracy and efficiency of numerical approximations for high-dimensional space FDEs.  Previous studies focusing on FDM and FEM methods have provided various accurate and mathematical convergent schemes for solving space FDE models \cite{Bu2015,Li2015,Zhao2015}, but the prohibitive computational cost challenges long-range time computation of the spatiotemporal FDEs in large spatial domains.

Various numerical technologies, in addition to analytical and approximate solutions \cite{Mainardi2007,Jiang2012}, have been employed to overcome the computational burden of the FDEs. Firstly, several approaches to reduce the computational cost of long-time range predictions for the temporal FDEs, have been considered \cite{Diethelm2006,Zeng2016}. For example, Yang et al. \cite{Yang2011} proposed a matrix transfer technique to reduce the computational cost of solving the two-dimensional spatiotemporal FDE.  Fu et al. \cite{Fu2013} introduced a Laplace transform technique combined with a boundary-only meshless particle method to effectively simulate the long time-history of fractional diffusion systems.  Sun et al. \cite{Sun2013} presented a semi-discrete FEM method to overcome the critical long-range-time computation problem for a class of time-fractional diffusion equations.  Secondly, extensive efforts have been made to reduce the computational cost caused by fractional spatial derivative terms in spatial FDE models.  For example, Feng et al. \cite{Feng2016} introduced a fast accurate iterative method for the Crank-Nicolson scheme which requires storage memory of $O(N)$ and the computational cost of $O(N\log N)$ for space-fractional diffusion equations, in which $N$ is the partition number in the space direction.  Wang and Basu \cite{Wang2012} developed a fast implicit finite difference discretization of two-dimensional space-fractional diffusion equations by decomposing the coefficient matrix into a combination of sparse and structured dense matrices.  Zhao et al. \cite{Zhao2016} presented a fast solver for space fractional differential equations using hierarchical matrices and a geometric multigrid method.  Pang and Sun \cite{Pang2016} introduced a fast numerical contour integral method for FDEs.  Zhang and Sun \cite{Zhang2015} applied the exponential quadrature rule for the coefficient matrix and obtained a fast numerical solution for FDEs.  Some other schemes, such as the FEM \cite{Du2015} and finite volume method \cite{Jia2016}, have also been proposed to solve the space FDEs by reducing the computational cost from $O(N^2)$ or $O(N^3)$ into $O(N\log N)$.  The above numerical schemes usually focus on spatial FDE models in regular domains without adopting fractional Laplacian operators.

Vector FDEs are needed to capture real-world diffusion in multi-dimensional media.  An extension of FDEs from one dimension to two or three dimensions is not as straightforward as that for the second-order diffusion equation, due to the global correlation nature of anomalous transport \cite{Meerschaert1999,Chen2016}.  Meerschaert et al. \cite{Meerschaert2004} proposed a vector fractional derivative as a generalization of the Laplacian term, which satisfies the global correlation property and can efficiently describe anomalous transport in complex media with preferential flow path.  Therefore, we adopt their definition and focus on the development of numerical methods for solving spatiotemporal FDEs containing the vector fractional derivative operator \cite{Du2015,Pang2015}.

In addition, real-world anomalous diffusion can cover a large domain, challenging classical grid-based numerical solvers such as the FDM.  The computational cost of the FEM on mesh generation and assembling of the asymmetric stiffness matrix also increases dramatically \cite{Liu2015}.  The radial basis function (RBF) collocation technique, which is a global meshless method, was recognized as an efficient method to approximate the differential equations in large and irregular domains \cite{LiuQ2015,Kansa2000,Chenw2014,Atluri2000}.  The RBF method has also proved to be a meshless merit method with the advantages of high accuracy and ease of implementation in discretizing the space-fractional derivative \cite{Pang2015,LiuQ2015,Dehghan2016}.

This study aims at providing an efficient method which can solve the spatiotemporal FDE in large, two-dimensional spatial domains, and provide a long-time range prediction of anomalous pollutant transport. We propose a semi-discrete Kansa method to solve the spatiotemporal FDE in large domains. Firstly, the Kansa method is employed for spatial discretization of the FDE in large spatial domains, and then the Gauss-Jacobi quadrature rule is used to accurately calculate the fractional derivative of the matrix obtained by the Kansa method. The computational cost of the spatiotemporal FDE model is then reduced since we only need to calculate the time fractional ordinary differential equation (ODE) in a matrix form. Finally, an exact solution for the ODE is obtained using the analytical technique and property of the Mittag-Leffler function, to approximate the time fractional derivative term in the spatiotemporal FDE model achieving a fast and accurate solution.

\section{Spatiotemporal FDE for anomalous transport in heterogeneous anisotropic media}
In the spatiotemporal FDE model, the time fractional derivative is used to describe the memory effect of solute dynamics.  The space derivative describes the spatial dependency or non-locality of solute particle displacement. Hereby, we consider the following two-dimensional spatiotemporal FDE
\begin{eqnarray}
\begin{array}{c}
\displaystyle{\frac{d^\alpha u(x,y,t)}{d t^\alpha} =-\textbf{V}(x,y) \cdot \nabla u(x,y,t)+k(x,y) D_M^\mathbf{\beta} u(x,y,t)+ f(x,y),}\\
(x,y)\in \Omega,\,\,\displaystyle{ 0<\alpha \leq 1,\,\,1<\mathbf{\beta} \leq 2,} \\
\displaystyle{u(x,y,0)=u_0(x,y),\,\, u(x,y,t)|_{\partial \Omega}=u_b(x,y,t),}
\end{array}
\label{eq0}
\end{eqnarray}
in which $u(x,y,t)$ represents the solute concentration, $k(x,y)$ is the diffusion coefficient, $\textbf{V}(x,y)$ is flow velocity, and $f$ denotes a source/sink term. Here the time fractional derivative with $\alpha$-order is defined in Caputo basis
\begin{eqnarray}
\begin{array}{c}
\displaystyle{\frac{d^\alpha u(t)}{d t^\alpha} =\frac{1}{\Gamma(1-\alpha)}\int_0^t \frac{u'(\tau)}{(t-\tau)^\alpha}d \tau,\,\,0<\alpha \leq 1.}
\end{array}
\label{eq1}
\end{eqnarray}

The generalized fractional derivative $D_M^\mathbf{\beta} (*)$ can be written as:
\begin{eqnarray}
\displaystyle{D_M^\beta u(x,y,t)=}\left\{
\begin{array}{c}
\displaystyle{\int_0^{2 \pi} D_\theta^{\beta(\theta)} u(x,y,t) m(\theta) d \theta\,\, (Continuous\,\,case),}\\
\displaystyle{\sum_i D_{\theta_i}^{\beta(\theta_i)} u(x,y,t) m(\theta_i) \,\, (Discrete\,\,case),}
\end{array}\right.\;
\label{eq1b}
\end{eqnarray}
where the weight function $m(\theta)$ or $m(\theta_i)$ (also called mixing measure) represents the probability of diffusion along the radial direction $\theta$ or $\theta_i$ \cite{Pang2015}. In the discrete case (i.e., with discontinuous angles $\theta_i$), the fractional derivative $\textbf{D}_{\theta_i}^{\beta_i}$ represents the vector Gr\"{u}nwald formula for fractional derivatives in different differentiation directions $\theta_i$. Here the order of the fractional derivative will be direction dependent, to capture the influence of multi-dimensional media with anisotropic heterogeneity on solute transport. The fractional derivative along the $\theta$ direction is defined as \cite{Roop2006,Meerschaert2004}:

\begin{eqnarray}
\begin{array}{c}
\displaystyle{\textbf{D}_{\theta}^\beta u=I_\theta^{2-\beta}D_\theta^2 u,}
\end{array}
\label{eq2}
\end{eqnarray}
where the symbol $``I"$ denotes the fractional integral. The above vector operator can be simplified for specific angles $\theta$ in an $x/y$ coordinate system; for example, when $\theta=0$, $\textbf{D}_{\theta=0}^\beta u$ reduces to $d^\beta/d x^\beta$; $\theta=\pi/2$, $\textbf{D}_{\theta=\pi/2}^\beta u$
reduces to $d^\beta/d y^\beta$. The corresponding fractional directional integral is given by
\begin{eqnarray}
\begin{array}{c}
\displaystyle{I_\theta^\alpha u(x,y)=\frac{1}{\Gamma (\alpha)}\int_0^{d(x,y,\theta)} \zeta^{\alpha-1} u(x-\zeta \cos \theta, y-\zeta \sin \theta) d \zeta, \; \alpha\in (0, 1),}
\end{array}
\label{eq2b}
\end{eqnarray}
where $d(x, y, \theta)$ is the distance from the internal node $(x, y)$ to the boundary $\partial \Omega$
along the direction $( -\cos\theta, -\sin \theta)$.

\section{Algorithm development} \label{sec:2}
\subsection{Preliminary knowledge}
\textbf{Mittag-Leffler function}

In our numerical scheme, the following one-parameter Mittag-Leffler function will be used:
\begin{eqnarray}
\begin{array}{c}
\displaystyle{E_\alpha (z)=\sum_{n=1}^\infty \frac{z^n}{\Gamma(\alpha n+1)},\,\, \alpha>0,\,\,z\in C}
\end{array}
\label{eq7}
\end{eqnarray}
where $\Gamma(z)$ represents the Gamma function which can be expressed as $\Gamma(z)=\int_0^\infty e^{-p} p^{z-1}dp$. The fractional derivative (Caputo type) of Mittag-Leffler function has the following properties \cite{Samko1993,Sun2013}:
\begin{eqnarray}
\left\{
\begin{array}{l}
\displaystyle{\frac{d^\alpha }{d t^\alpha} E_\alpha (\lambda t^\alpha)= \lambda E_\alpha (\lambda t^\alpha),} \\
\displaystyle{\frac{d^\alpha C}{d t^\alpha}  = 0,\,\, C\,\, } \mbox{is a constant}.
\end{array}\right.\;
\label{eq8}
\end{eqnarray}

\textbf{Kansa method}

In the Kansa method, the solute concentration in Eq. (\ref{eq0}) is approximated by the weighted summation of the RBF $\phi(r)$ \cite{Kansa1990a,Kansa1990b,Chen2010}:
\begin{eqnarray}
\begin{array}{c}
\displaystyle{u(x,y,t)\approx \sum_{i=1}^{M+N}\lambda_i(t)\phi (r_i) \,\, ,}
\end{array}
\label{eq5}
\end{eqnarray}
where $r_j=\sqrt{(x-x_i)^2+(y-y_i)^2}$; $(x_i, y_i)$ is the source point in the computational domain; $M$ and $N$ are the number of domain points and boundary points, respectively; and $\lambda_i (t)$ is the unknown expansion coefficient. There are several kinds of RBFs employed in the Kansa method, and here we use a popular and efficient one called Hardy's multiquadrics (MQ)
\begin{eqnarray}
\begin{array}{c}
\displaystyle{\phi (r_i)=\sqrt{r_i^2+C^2} \,\, ,}
\end{array}
\label{eq6}
\end{eqnarray}
where the shape parameter $C$ ($C>0$) can be determined using algorithms developed in the literature \cite{Hardy1971,Fasshauer}.

\subsection{Numerical scheme}
For illustration purposes, here we consider the space fractional derivative in orthogonal coordinates $x$ and $y$, which can be simply written as $D_{\theta=0}^\beta$ and $D_{\theta=\pi/2}^\beta$. The resultant governing equation is
\begin{eqnarray}
\left\{
\begin{array}{c}
\displaystyle{\frac{d^\alpha u(x,y,t)}{d t^\alpha} =-\mathbf{V}(x,y) \cdot \nabla u(x,y,t)+k_x(x,y) \frac{\partial^{\beta_x}u(x,y,t)}{\partial x^{\beta_x}}}\\
\displaystyle{+k_y(x,y) \frac{\partial^{\beta_y}u(x,y,t)}{\partial y^{\beta_y}}
+ f(x,y),\,\,in\,\, \Omega \times (0, T),}\\
\displaystyle{u(x,y,0)=u_0(x,y),\,\, u(x,y,t)|_{\partial \Omega}=u_b(x,y,t),}
\end{array}
\right.\;
\label{eq3}
\end{eqnarray}
where $0<\alpha \leq 1,\,\,1<\beta_x,\,\,\beta_y\leq 2$, $k_x(x,y)$ and $k_y(x,y)$ are dispersion coefficients, $f(x,y,t)$ denotes the source/sink term, and $V(x,y)$ is the water flow velocity with components $V_x(x,y)$ and $V_y(x,y)$. The advection term can be written as
\begin{eqnarray}
\begin{array}{c}
\displaystyle{\mathbf{V}(x,y) \cdot \nabla u(x,y,t)=V_x(x,y) \frac{\partial u(x,y,t)}{\partial x}+V_y(x,y) \frac{\partial u(x,y,t)}{\partial y}}.
\end{array}
\label{eq4}
\end{eqnarray}

Using the interpolation formula Eq. (\ref{eq5}), a set of semi-discrete discretization of Eq. (\ref{eq3}) which satisfies the governing equation, and initial and boundary conditions at $M+N$ collocation points can be written as
\begin{eqnarray}
\left\{
\begin{array}{c}
\displaystyle{\sum_{j=1}^{M+N} \phi_{ij} \frac{d^\alpha \lambda_j(t)}{dt^\alpha}=-V_x(x_i,y_i) \sum_{j=1}^{M+N} (\frac{\partial \phi}{\partial x}) _{i j} \lambda_j(t)-V_y(x_i,y_i) \sum_{j=1}^{M+N} (\frac{\partial \phi}{\partial y}) _{i j} \lambda_j(t)}\\
\displaystyle{+k_{x}(x_i, y_i) \sum_{j=1}^{M+N} (\frac{\partial^{\beta_x}\phi}{\partial x^{\beta_x}})_{ij}\lambda_j(t)+k_{y}(x_i,y_i) \sum_{j=1}^{M+N} (\frac{\partial^{\beta_y}\phi}{\partial y^{\beta_y}})_{ij} \lambda_j(t)
+ f_i,\,\,i=1,2,...M,}\\
\displaystyle{\sum_{j=1}^{M+N} \phi_{ij}\lambda_j (t)=(u_b)_i,\,\,i=M+1,M+2,...M+N,}\\
\displaystyle{\sum_{j=1}^{M+N} \phi_{ij}\lambda_j (t=0)=(u_0)_i,\,\,i=1,2,...M+N.}
\end{array}
\right.\;
\label{eq10}
\end{eqnarray}
Here, the subscript in $\phi_{i j}$ and $f_i$ represents the sequence of points for evaluating the function values. Discretization of Eq. (\ref{eq10}) leads to the following matrix form
\begin{eqnarray}
\left\{
\begin{array}{c}
\displaystyle{\mathbf{\Phi}_d \frac{d^\alpha \mathbf{\lambda}}{dt^\alpha}=[- \mathbf{V}_x\frac{\partial \mathbf{\Phi}_d}{\partial x}- \mathbf{V}_y \frac{\partial \mathbf{\Phi}_d}{\partial y}
+\mathbf{K}_{x} (\frac{\partial^{\beta_x} \mathbf{\Phi}}{\partial x^{\beta_x}})+\mathbf{K}_{y} (\frac{\partial^{\beta_y} \mathbf{\Phi}}{\partial y^{\beta_y}})] \mathbf{\lambda}
+ \mathbf{F},}\\
\displaystyle{0=-\mathbf{\Phi}_b \mathbf{\lambda}+\mathbf{u}_b,}\\
\displaystyle{\mathbf{\Phi} \mathbf{\lambda}_{(t=0)}=\mathbf{u}_0,}
\end{array}
\right.\;
\label{eq11}
\end{eqnarray}
where $\mathbf{\Phi}$ is a $M\times (M+N)$ coefficient matrix; $\mathbf{\Phi}_d=\{\phi_{ij}\}$ for $i=1,2,...,M$ and $j=1,2,...,M+N$; $\mathbf{\Phi}_b=\{\phi_{ij}\}$ for $i=M+1,M+2,...,M+N$ and $j=1,2,...,M+N$;
$\mathbf{\lambda}^{(n)}=\{\lambda_j [(n-1)\Delta t]\}$ for $j=1,2,...,M+N$; and $\mathbf{V}_x$, $\mathbf{V}_y$, $\mathbf{K}_x$, $\mathbf{K}_y$ and $\mathbf{F}$ are $M\times (M+N)$
 matrices with nonzero entries $\{V_{ii}\}$, $\{k_{xii}\}$ and $\{k_{yii}\}$.

In the next step, we evaluate the spatial derivative of the RBFs including $(\partial^{\beta_x}\phi/\partial x^{\beta_x})_{ij}$,
$(\partial^{\beta_y}\phi/\partial y^{\beta_y})_{ij}$ and the advection term. We first consider the term $(\partial^{\beta_x}\phi/\partial x^{\beta_x})_{ij}$ as an example:
\begin{eqnarray}
\begin{array}{c}
\displaystyle{(\frac{\partial^{\beta_x}\phi}{\partial x^{\beta_x}})_{ij}=\frac{1}{\Gamma(2-\beta_x)}\int_0^{x_i} (x_i-\varsigma)^{1-\beta_x}\frac{\partial^2 \phi(\sqrt{(x-x_j)^2+(y_i-y_j)^2})}{\partial x^2}d \varsigma}.
\end{array}
\label{eq12}
\end{eqnarray}
The final approximation of the above derivative term can be written as follows using variable substitution and the Gauss-Jacobi quadrature method \cite{Pang2015}
\begin{eqnarray}
\begin{array}{c}
\displaystyle{(\frac{\partial^{\beta_x}\phi}{\partial x^{\beta_x}})_{ij}=\frac{(x_i/2)^{2-\beta_x}}{\Gamma(2-\beta)} \sum_k \omega_k \chi(\xi_k)},
\end{array}
\label{eq13}
\end{eqnarray}
where $\xi_k$ and $\omega_k$ are Gauss-Jacobi quadrature points and weights, respectively; and
\begin{eqnarray}
\begin{array}{c}
\displaystyle{\chi(\xi_k)=\frac{\partial^2 \phi(\sqrt{(x-x_j)^2+(y_i-y_j)^2})}{\partial x^2}|_{x=x_i-x_i(1+\xi_k)/2}}.
\end{array}
\label{eq14}
\end{eqnarray}
Details of this approximation can be found in Ref.\cite{Pang2015}. Approximation of the spatial derivative along the $y$ direction and the advection term can be achieved following the same argument.

By substituting the above approximation into Eq. (\ref{eq11}), and assuming
\begin{eqnarray}
\mathbf{C}=\Phi=\left[
\begin{array}{c}
     \mathbf{\Phi}_d\\
     0
\end{array}\right]_{(M+N)\times (M+N)},\;
\label{eq15}
\end{eqnarray}

 \begin{eqnarray}
\mathbf{B}=\left[
\begin{array}{c}
     \mathbf{V}_x \frac{\partial \mathbf{\Phi}_d}{\partial x}+\mathbf{V}_y \frac{\partial \mathbf{\Phi}_d}{\partial y}
-\mathbf{K}_{x} (\frac{\partial^{\beta_x} \mathbf{\Phi}}{\partial x^{\beta_x}})-\mathbf{K}_{y} (\frac{\partial^{\beta_y} \mathbf{\Phi}}{\partial y^{\beta_y}})\\
     \Phi_b
\end{array}\right]_{(M+N)\times (M+N)},\;
\label{eq16}
\end{eqnarray}
 and
 \begin{eqnarray}
\mathbf{\overline{F}}=\left[
\begin{array}{c}
     \mathbf{F}\\
     \Phi_b
\end{array}\right]_{(M+N)\times (M+N)},\;
\label{eq17}
\end{eqnarray}
we obtain the following equations:
\begin{eqnarray}
\begin{array}{c}
\displaystyle{\mathbf{C}\frac{\partial^{\alpha}\lambda_t}{\partial t^{\alpha}}=-\mathbf{B} \lambda_t+\mathbf{\overline{F}}}.
\end{array}
\label{eq18}
\end{eqnarray}
Let $\overline{\lambda_t}=\lambda_t-\mathbf{B}^{-1} \mathbf{\overline{F}}$, one can obtain the matrix form of Eq. (\ref{eq11}):
%\begin{eqnarray}
%\begin{array}{c}
%\displaystyle{\mathbf{C}\frac{\partial^{\alpha}\overline{\lambda_t}}{\partial t^{\alpha}}=-\mathbf{B} \overline{\lambda_t}}.
%\end{array}
%\label{eq19}
%\end{eqnarray}
\begin{eqnarray}
\left\{
\begin{array}{c}
\displaystyle{\mathbf{C}\frac{\partial^{\alpha}\overline{\lambda_t}}{\partial t^{\alpha}}=-\mathbf{B} \overline{\lambda_t}}\\
\displaystyle{\overline{\lambda_t}_{(t=0)}=\overline{\mathbf{u}}_0,}
\end{array}
\right.\;
\label{eq20}
\end{eqnarray}
where $\overline{\mathbf{u}}_0=\Phi^{-1}\mathbf{u}_0-\mathbf{B}^{-1} \mathbf{\overline{F}}$.

The exact solution of the fractional relaxation equation (\ref{eq20})
can be expressed as \cite{Sun2013,Mainardi1996}:
\begin{eqnarray}
\begin{array}{c}
\displaystyle{\overline{\lambda_t}= \overline{\mathbf{u}}_0 E_\alpha (-\textbf{M}t^\alpha)}
\end{array}
\label{eq21}
\end{eqnarray}
where $\textbf{M}=\textbf{C}^{-1}\textbf{B}$, and $E_\gamma (*)$ is the
Mittag-Leffler function which can be accurately evaluated with
an existing Matlab code \cite{Podlubny2009b}. In our tests, the value of this
function was accurate within an order of $10^{-12}$.

Next, we decompose the Mittag-Leffler function in (\ref{eq21})  as:
\begin{eqnarray}
\begin{array}{c}
\displaystyle{E_\alpha (-\textbf{M} t^\alpha)=\textbf{S} \Lambda_t
\textbf{S}^{-1}}
\end{array}
\label{eq22}
\end{eqnarray}
where $\textbf{S}$ is the modal matrix formed by the eigenvectors of
$-\textbf{M}$. $\Lambda_t$ is a diagonal matrix
whose $i$-th diagonal entry is $E_\gamma (\Lambda_i t^\gamma)$,
and $\Lambda_i$ is the $i$-th eigenvalue of $-\textbf{M}$.
Inserting (\ref{eq22}) into (\ref{eq21}), we get
\begin{eqnarray}
\begin{array}{c}
\displaystyle{\tilde{\mathbf{\lambda}}_t= \textbf{S} \Lambda_t
\textbf{S}^{-1} \overline{\mathbf{u}}_0}.
\end{array}
\label{eq23}
\end{eqnarray}

After the above manipulations, the initial-boundary value FDE model
 (\ref{eq3}) is reduced to an ordinary differential system through RBF discretization. The reduced system can be solved
analytically in terms of the Mittag-Leffler function. Finally, an accurate numerical approximation of Eq. (\ref{eq3}) can be achieved via the RBF interpolation formula (\ref{eq5}).

For the general fractional derivative $D_M^\beta \: (\beta\in (1,2])$ along an arbitrary direction $\theta$, the semi-discrete discretization of Eq. (\ref{eq10}) can be rewritten as
\begin{eqnarray}
\left\{
\begin{array}{c}
\displaystyle{\sum_{j=1}^{M+N} \phi_{ij} \frac{d^\alpha \lambda_j(t)}{dt^\alpha}=-\sum_{j=1}^{M+N} V_i(x,y)(\frac{\partial \phi}{\partial x}+\frac{\partial \phi}{\partial y}) _{i j} \lambda_j(t)}\\
\displaystyle{+\sum_{j=1}^{M+N} k_{i}(x,y) (D_M^\beta \phi)_{ij}\lambda_j(t)
+ f_i,\,\,i=1,2,...M,}\\
\displaystyle{\sum_{j=1}^{M+N} \phi_{ij}\lambda_j (t)=(u_b)_i,\,\,i=M+1,M+2,...M+N,}\\
\displaystyle{\sum_{j=1}^{M+N} \phi_{ij}\lambda_j (t=0)=(u_0)_i,\,\,i=1,2,...M+N,}
\end{array}
\right.\;
\label{eq24}
\end{eqnarray}
where the general fractional derivative of the basis function $\phi$ is expressed as:
\begin{eqnarray}
\displaystyle{(D_M^\beta \phi)_{ij}=}\left\{
\begin{array}{c}
\displaystyle{\Sigma_l (D_{\theta_l}^{\beta} \phi)_{ij} \omega_l m(\theta_l)\,\, (Continuous\,\,case),}\\
\displaystyle{\Sigma_l (D_{\theta_l}^{\beta} \phi)_{ij} m_l \,\, (Discrete\,\,case).}
\end{array}\right.\;
\label{eq25}
\end{eqnarray}
Here $\theta_l$ and $\omega_k$ are the quadrature points and weights of the Gauss-Jacobi rule which approximates the outer integral with respect to direction $\theta$.
The fractional directional derivative of a RBF $\phi$, using the Kansa method (\ref{eq5}), can be written as
\begin{eqnarray}
\begin{array}{c}
\displaystyle{(D_{\theta}^{\beta} \phi)_{ij}=\frac{(d_{ii}(\theta)/2)^{2-\beta}}{\Gamma(2-\beta)} \sum_k \omega_k \chi(\xi_k)},
\end{array}
\label{eq26}
\end{eqnarray}
where $d_{ij}$ is the distance from an internal node to the domain boundary \cite{Pang2015}, $\omega_k$ is the Gauss-Jacobi quadrature point, and $\chi(\xi)$ is defined by
\begin{eqnarray}
\begin{array}{c}
\displaystyle{\chi(\xi)=\frac{1}{(r_{ij}^2(\theta)+C^2)^{1/2}}-\frac{A^2}{(r_{ij}^2(\theta)+C^2)^{3/2}}}.
\end{array}
\label{eq27}
\end{eqnarray}
The coefficients $A$ and $r_{ij}(\theta)$ are expressed as below,
\begin{eqnarray}
\begin{array}{c}
\displaystyle{A=(x_i-x_j) \cos\theta+(y_i-y_j) \sin\theta-\frac{r_{ij}(\theta)}{2}(1+\xi)},
\end{array}
\label{eq29}
\end{eqnarray}

\begin{eqnarray}
\begin{array}{c}
\displaystyle{r_{ij}(\theta)=\sqrt{[x_i- \cos \theta d_{ii}(\theta)(1+\xi)/2-x_j]^2+[y_i- \sin \theta d_{ii}(\theta)(1+\xi)/2-y_j]^2}}.
\end{array}
\label{eq29b}
\end{eqnarray}

\section{Numerical examples, validation, and analysis}
\label{sec:5}
\subsection{Example 1: One-dimensional spatiotemporal FDE}
\begin{figure}
  \includegraphics[width=0.8\textwidth]{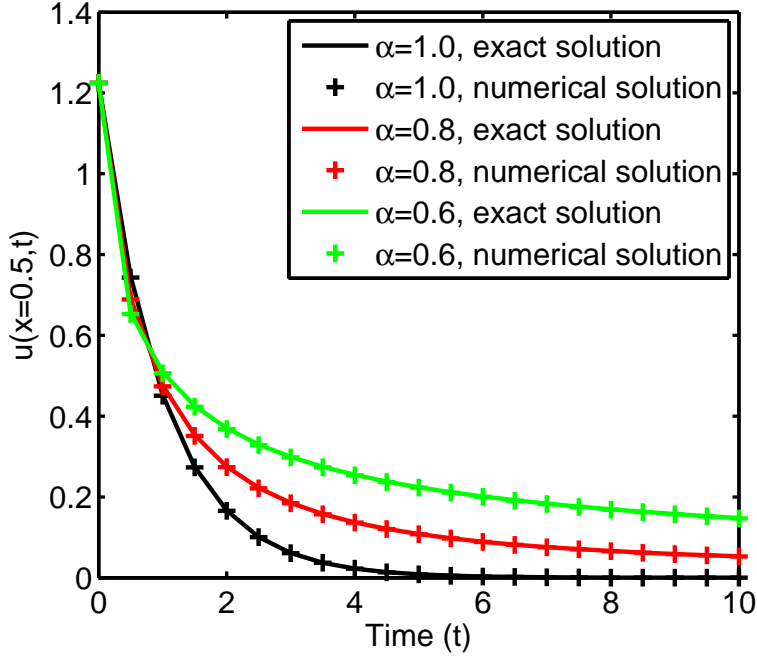}
\caption{A comparison of exact and numerical solutions of the spatiotemporal FDE model at $x=0.5$ with the order of the spatial derivative $\beta=1.6$. Here, the spatial nodes are
uniformly distributed with an interval $\Delta x=0.05$, and the shape parameter $C=0.1$.}
\label{fig:1}       % Give a unique label
\end{figure}

For validation purposes, we first consider a one-dimensional simplification of the spatiotemporal FDE (\ref{eq0}):
\begin{eqnarray}
\left\{
\begin{array}{c}
\displaystyle{\frac{\partial^\alpha u(x,t)}{\partial t^\alpha} =-V(x)\frac{\partial u(x,t)}{\partial x}+K \frac{\partial^\beta u(x,t)}{\partial x^\beta},}\\
\displaystyle{x\in (0, 1),\,\,t\in (0, 10],\,\,0<\alpha \leq 1,\,\, 1<\beta\leq 2,}\\
\displaystyle{u(0,t) =E_\alpha (-t^\alpha),\,\,u(1,t)=\sqrt{2} E_\alpha (-t^\alpha),}\\
\displaystyle{u(x,0) =(x+1)^{1/2},\,\, x\in [0, 1],}
\end{array}
\right.\;
\label{eq30}
\end{eqnarray}
where $K=-\Gamma (3/2-\beta)(x+1)^\beta/2\Gamma (3/2)$, $V(x)=x+1$, and the exact solution can be written as $u(x,t)=(x+1)^{1/2} E_\alpha (-t^\alpha)$.
An observation of Fig. \ref{fig:1} shows that the proposed scheme agrees well with the analytical solution.

To further explore the accuracy and convergence of the numerical scheme, we evaluate the impact of the spatial interval (between the uniformly distributed nodes) on numerical solutions. Tab. \ref{tab1} lists the maximum absolute error (MAE) ($l_\infty$) and the convergence rate of the Kansa method. Here the convergence rate is calculated in terms of the ratio of MAE at different grid sizes. Results show that the Kansa method has a super-linear convergence rate, which is consistent with our previous work on space fractional diffusion model \cite{Pang2015}. Moreover, the relative errors presented in Tab. \ref{tab2} verify that the proposed scheme can produce accurate long-range time computation without causing a large computational burden, due to the analytical approach used for the temporal terms in the model.

\begin{table}
\centering
\caption{The maximum absolute error (ABE) for grids refined for Example 1 at time $t=10$, with the order of the time fractional derivative $\alpha=0.6$ and the order of the space fractional derivative $\beta=1.6$. In the numerical computation, the accuracy rate of the Mittag-Leffler function is $10^{-12}$, and the shape parameter $C=0.1$.}
\begin{tabular}{lll}
\hline\noalign{\smallskip}
$\Delta x$ & $l_\infty$ error & Error rate \\
\noalign{\smallskip}\hline\noalign{\smallskip}
$\Delta x=1/10$  & 0.00494104 &    \\
$\Delta x=1/20$  & 0.00140965 & $3.51>20/10$ \\
$\Delta x=1/25$  & 0.00082844 & $1.70>25/20$ \\
$\Delta x=1/50$  & 0.00029132 & $2.84>50/25$\\
\noalign{\smallskip}\hline
\end{tabular}
\label{tab1}
\end{table}
\begin{table}
\caption{The relative error (Error$=|(u_{exact}(L/2, t)-u(L/2, t))/u_{exact}(L/2, t)|$) and the CPU time of a PC for different total modeling times with the domain length of $L=1.0$, the node spacing of $\Delta x=1/20$, and the shape parameter $C=0.1$.}
% For LaTeX tables use2
\begin{tabular}{lllll}
\hline\noalign{\smallskip}
$Time\,\, (t) $ & $t=100$ &$t=1000$ &$t=10000$ & $t=100000$ \\
\noalign{\smallskip}\hline\noalign{\smallskip}
Relative error  & 0.02243720 & 0.02243706 & 0.02243702&  0.02243701 \\
\noalign{\smallskip}\hline\noalign{\smallskip}
CPU time & $0.043852 s$ & $0.051386 s$ & $0.052045 s$ & $0.066588 s$\\
\noalign{\smallskip}\hline
\end{tabular}
\label{tab2}
\end{table}

\subsection{Example 2: Two-dimensional spatiotemporal FDE}

\textbf{a. Rectangular domain}
The example problem is defined as:
\begin{eqnarray}
\left\{
\begin{array}{c}
\displaystyle{\frac{\partial^\alpha u(x,y,t)}{\partial t^\alpha} =K_x(x,y) \frac{\partial^{\beta_x} u(x,y,t)}{\partial x^{\beta_x}}+K_y(x,y) \frac{\partial^{\beta_y} u(x,y,t)}{\partial y^{\beta_y}},}\\
\displaystyle{(x,y)\in (0, 1)\times (0, 1),\,\,t\in (0, 10],\,\,0<\alpha \leq 1,\,\, 1<\beta_x,\,\,\beta_y\leq 2,}\\
\displaystyle{u(0,y,t) =(y+1)^{1/2} E_\alpha (-t^\alpha),\,\,u(1,y,t)=\sqrt{2} (y+1)^{1/2} E_\alpha (-t^\alpha),}\\
\displaystyle{u(x,0,t) =(x+1)^{1/2} E_\alpha (-t^\alpha),\,\,u(x,1,t)=\sqrt{2} (x+1)^{1/2} E_\alpha (-t^\alpha),}\\
\displaystyle{u(x,y,0) =(x+1)^{1/2} (y+1)^{1/2},\,\, (x,y)\in (0, 1)\times (0, 1).}
\end{array}
\right.\;
\label{eq31}
\end{eqnarray}
If $K_x=-2\Gamma (3/2-\beta_x)(x+1)^{\beta_x}/3\Gamma (3/2)$ and $K_y=-\Gamma (3/2-\beta_y)(y+1)^{\beta_y}/3\Gamma (3/2)$, then the analytical solution is $u(x,t)=(x+1)^{1/2} (y+1)^{1/2} E_\alpha (-t^\alpha)$.

\begin{figure}[htb]
\centering \subfigure{
\includegraphics[width=0.45\linewidth]{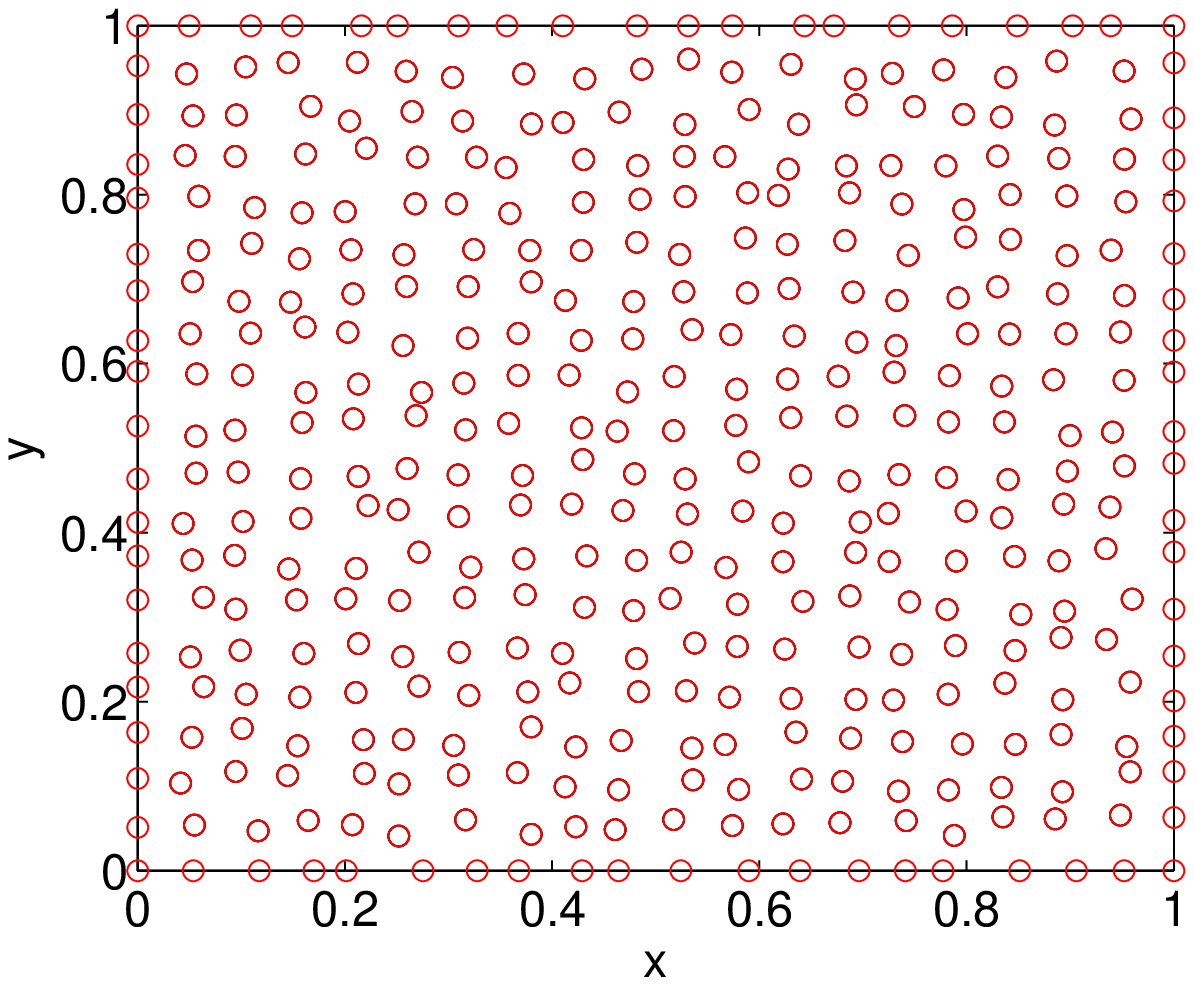}
\label{fig:forceBalanceComparison_a}} \subfigure{
\includegraphics[width=0.45\linewidth]{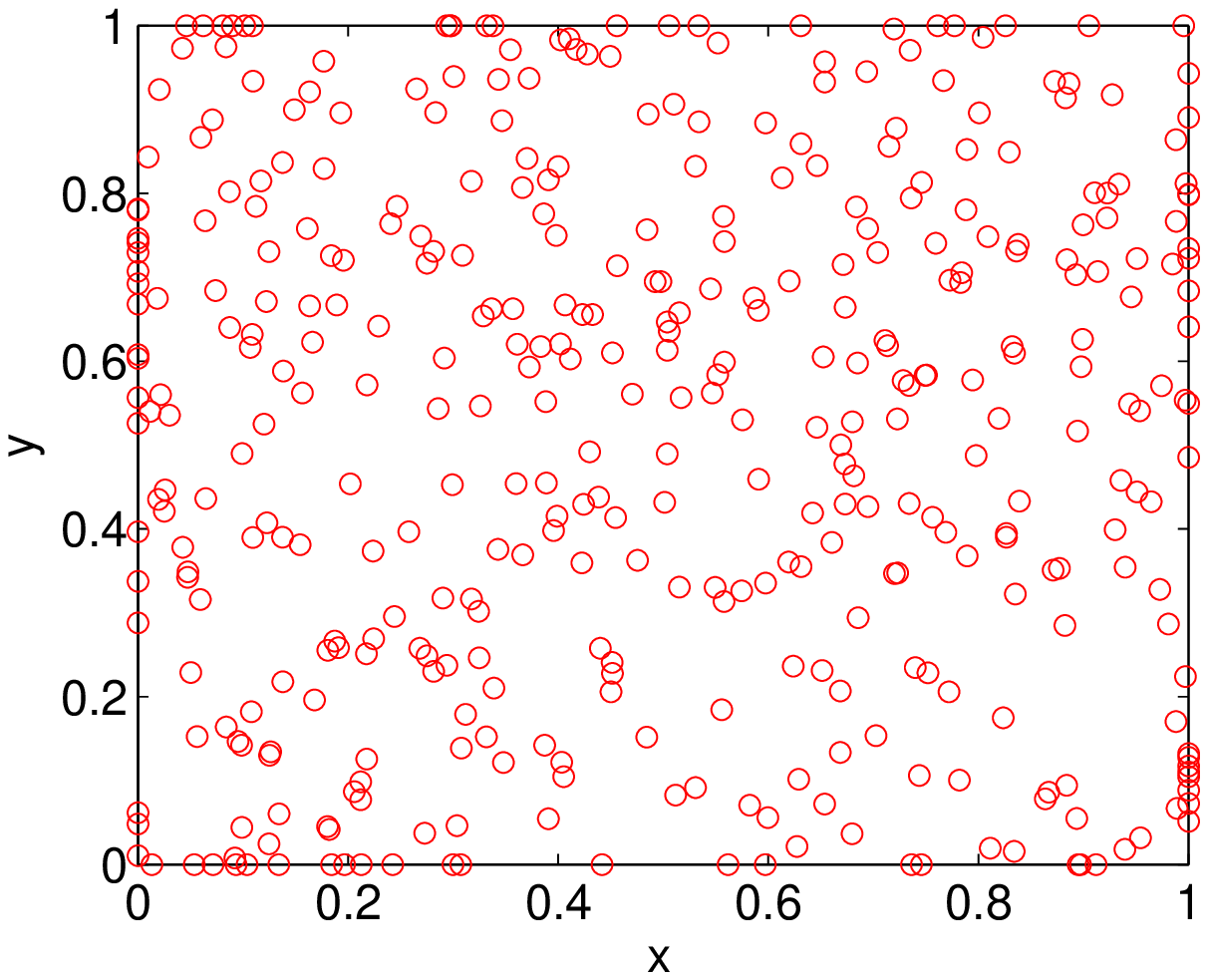}
\label{fig:forceBalanceComparison_b}} \caption{Configuration of irregular nodes in Example 2, type(I) (left) and type(II) (right).} \label{fig3}
\end{figure}

Because one of the advantages of the RBF method is that it allows a random distribution of collocation nodes, we can investigate the influence of the node distribution modes on the accuracy of numerical solutions. Here we consider two types of regional node-distribution modes: (I) Irregular modes with nodes derived by jiggling the regular nodes, i.e. assigning a perturbation on regular nodes along the $x$ and $y$ directions, as shown in Fig.\ref{fig3} (Left); (II) Irregular modes with nodes generated by the random vector of a uniform distribution, as shown in Fig. \ref{fig3} (Right).

\begin{figure}[htb]
\centering \subfigure{
\includegraphics[width=0.8\linewidth]{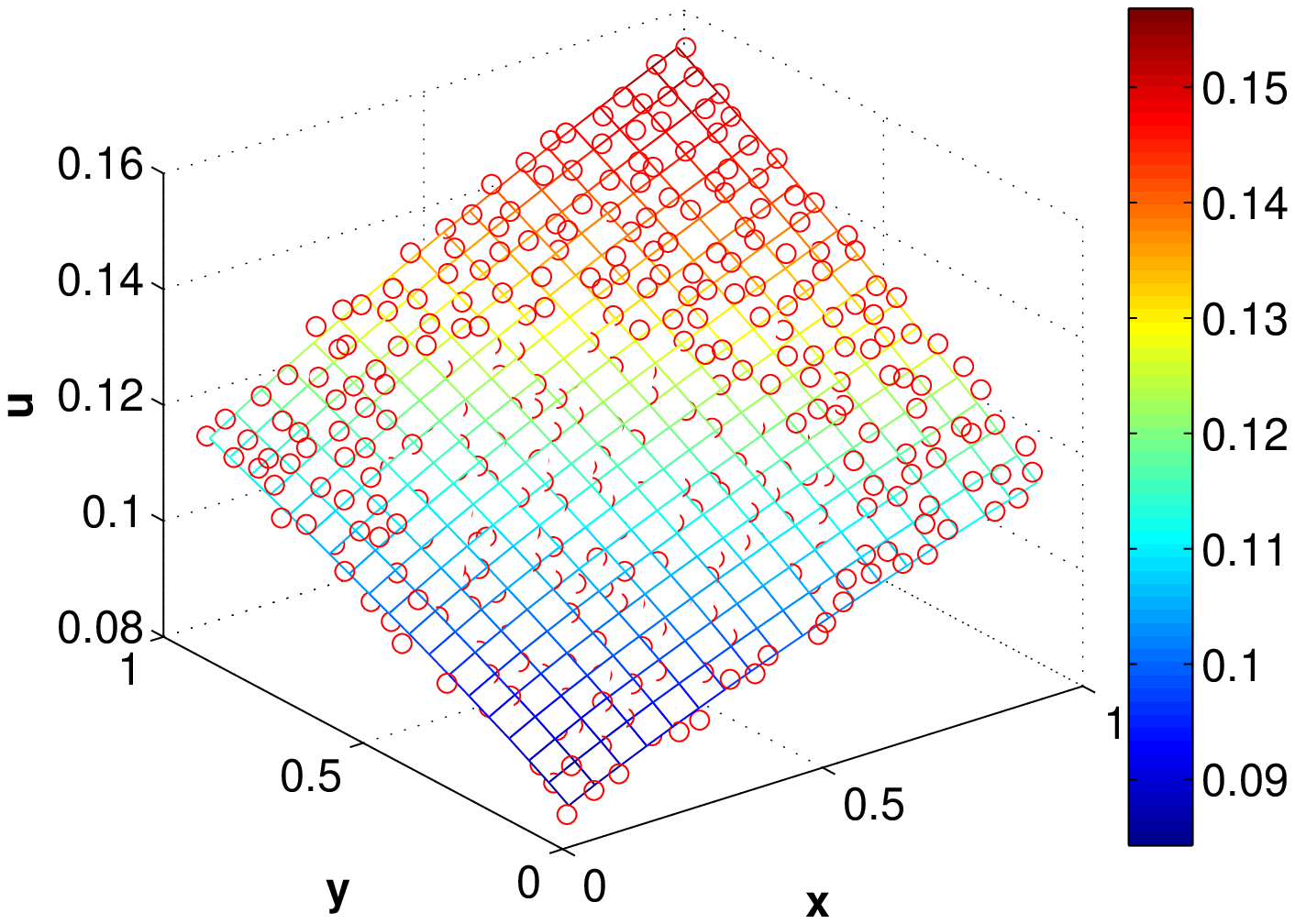}
\label{fig:forceBalanceComparison_a}} \subfigure{
\includegraphics[width=0.8\linewidth]{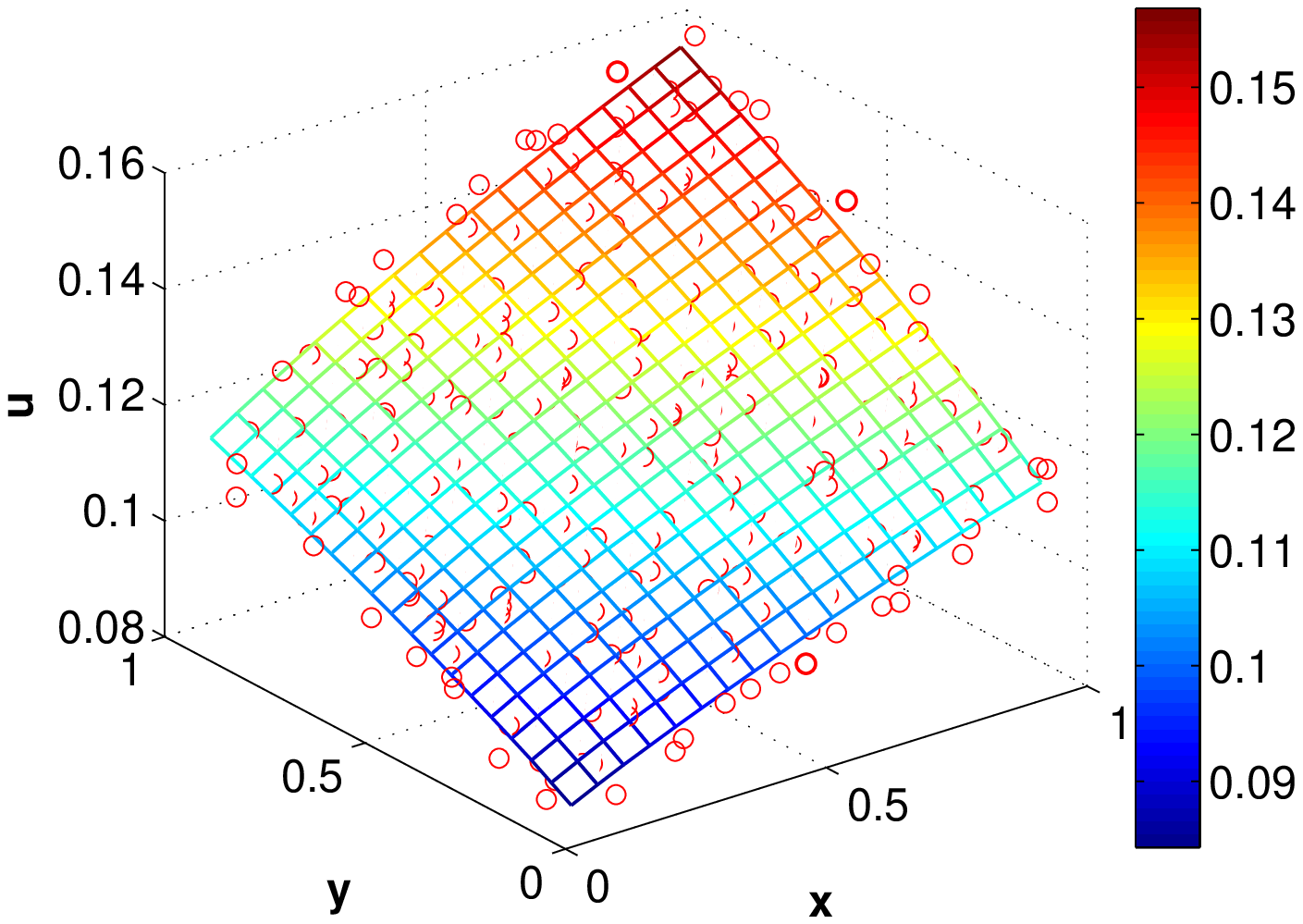}
\label{fig:forceBalanceComparison_b}} \caption{A comparison of the numerical results with irregular nodes type(I) and (II) and the analytical solutions for the spatiotemporal FDE (Eq. \ref{eq31})
at time $t=10$ with $\alpha=0.7$, $\beta_x=1.6$, and $\beta_y=1.8$. In the numerical scheme, $\Delta x=\Delta y=1/20$, and the MQ parameter $c=0.01$.} \label{fig5}
\end{figure}

A comparison of the numerical results using irregular nodes (i.e., types (I) and (II)) with the exact solutions is presented in Fig. \ref{fig5}.  Results show that both node-distribution modes offer relatively good numerical solutions for the presented equation. However, the randomly distributed nodes in type (II) may cause ill-condition problems in some cases.

\begin{figure}
  \includegraphics[width=0.8\textwidth]{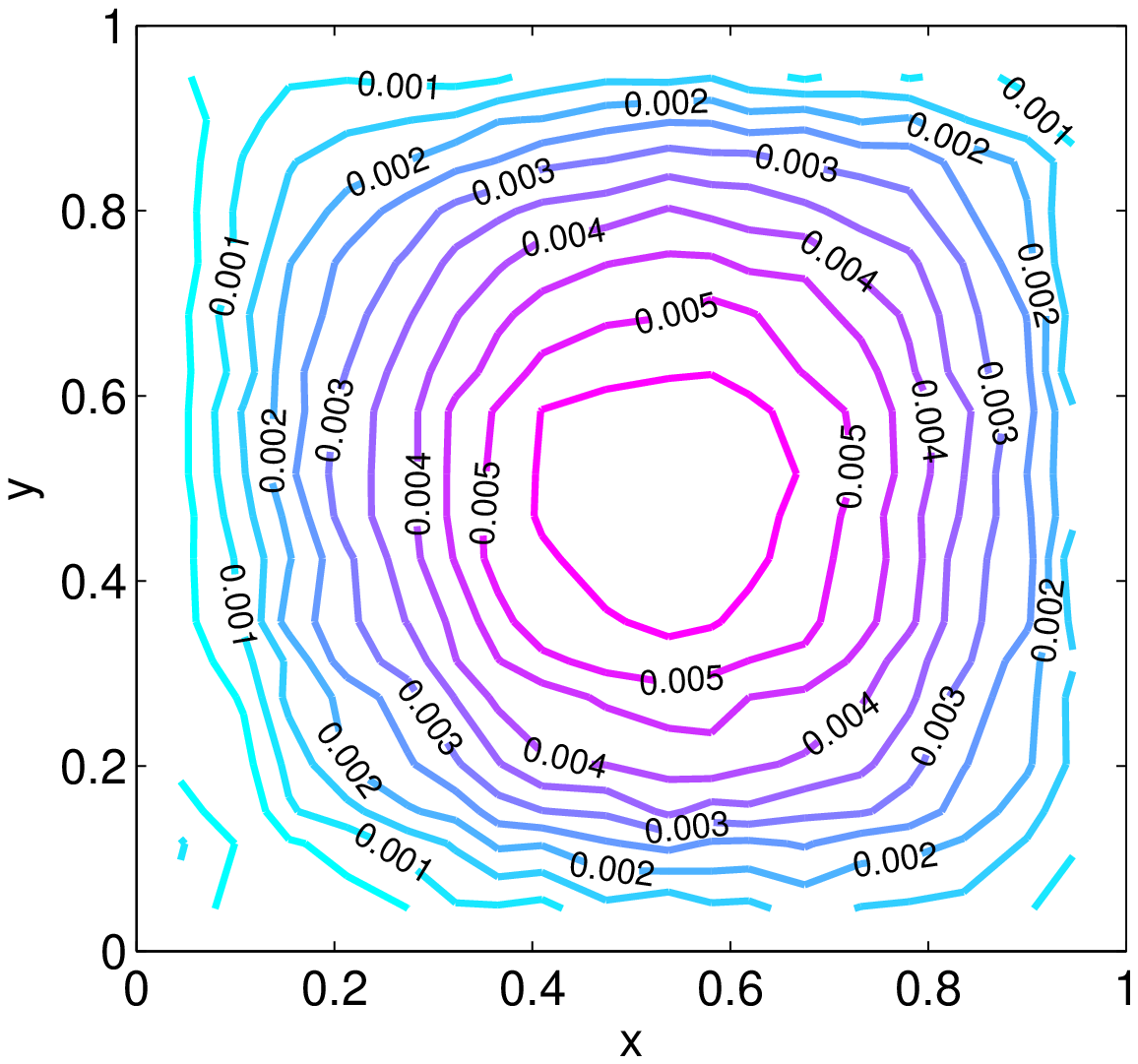}
\caption{Contour map of the relative error with irregular grids (type(I)) for the spatiotemporal FDE (Eq. \ref{eq31}) at time $t=10$ with $\alpha=0.7$, $\beta_x=1.6$, and $\beta_y=1.8$. In the numerical scheme $\Delta x=\Delta y=1/20$, and the MQ parameter $c=0.01$.}
\label{fig6}       % Give a unique label
\end{figure}

If we use irregular grids (type (I)), the contour map of the relative errors (Fig.\ref{fig6}) shows that relative errors increase with the distance between each node and the boundary domain, due to the first-type boundary condition employed in the present example.

\begin{table}
\centering
\caption{MAEs and relative errors ($max(|(u_{exact}(:, t)-u(:, t))/u_{exact}(:, t)|)$) for regular grids refined for Example 2 at $t=10$ with $\alpha=0.7$, $\beta_x=1.6$, and $\beta_y=1.8$. In the numerical scheme, the MQ parameter $c=0.01$.}
\begin{tabular}{llll}
\hline\noalign{\smallskip}
$\Delta =\Delta x=\Delta y$ & $l_{\infty}$ error & Error rate & Relative error  \\
\noalign{\smallskip}\hline\noalign{\smallskip}
$\Delta =1/10$  & 0.01003497 &   & 0.01793833 \\
$\Delta =1/15$  & 0.00469933 & $2.14>15/10$ & 0.00818898\\
$\Delta =1/20$  & 0.00281028 & $1.67>20/15$ & 0.00534286\\
$\Delta =1/25$  & 0.00278984 & $1.01<25/20$ & 0.00593439\\
\noalign{\smallskip}\hline
\end{tabular}
\label{tab3}
\end{table}

Tab. \ref{tab3} shows that the discretization of the spatial derivative terms using the Kansa method provides reasonable numerical results. For example, the MAE is less than $10^{-2}$, and the maximum relative error is less than $1\%$ for $\Delta x=\Delta y=1/15$.  Moreover, the MAE decreases with an increase of the collocation node number. However, Tab. \ref{tab3} also shows that the error rate of the Kansa method does not always decline faster than linear with an increase in the collocation node number. Instead, the error rate may decrease more slowly with an increasing collocation number in many cases \cite{Pang2015,Chenw2014}. Therefore we should consider both accuracy and computational cost in real-world applications when employing the Kansa method. ~\\

\textbf{b. Unit disk domain}

\begin{figure}[htb]
\centering
\subfigure{
\includegraphics[width=0.45\linewidth]{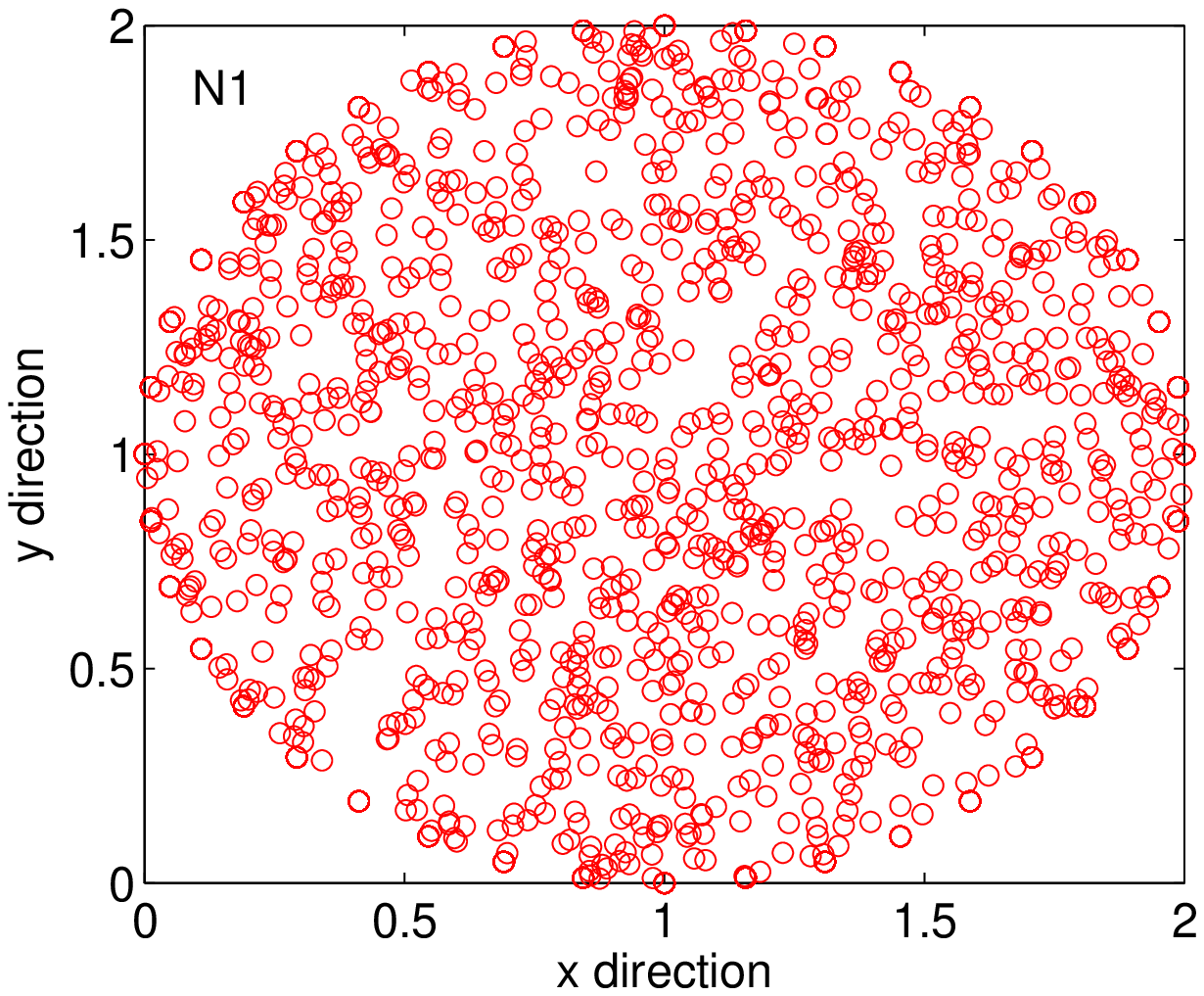}
\label{fig:forceBalanceComparison_a}}
\subfigure{
\includegraphics[width=0.45\linewidth]{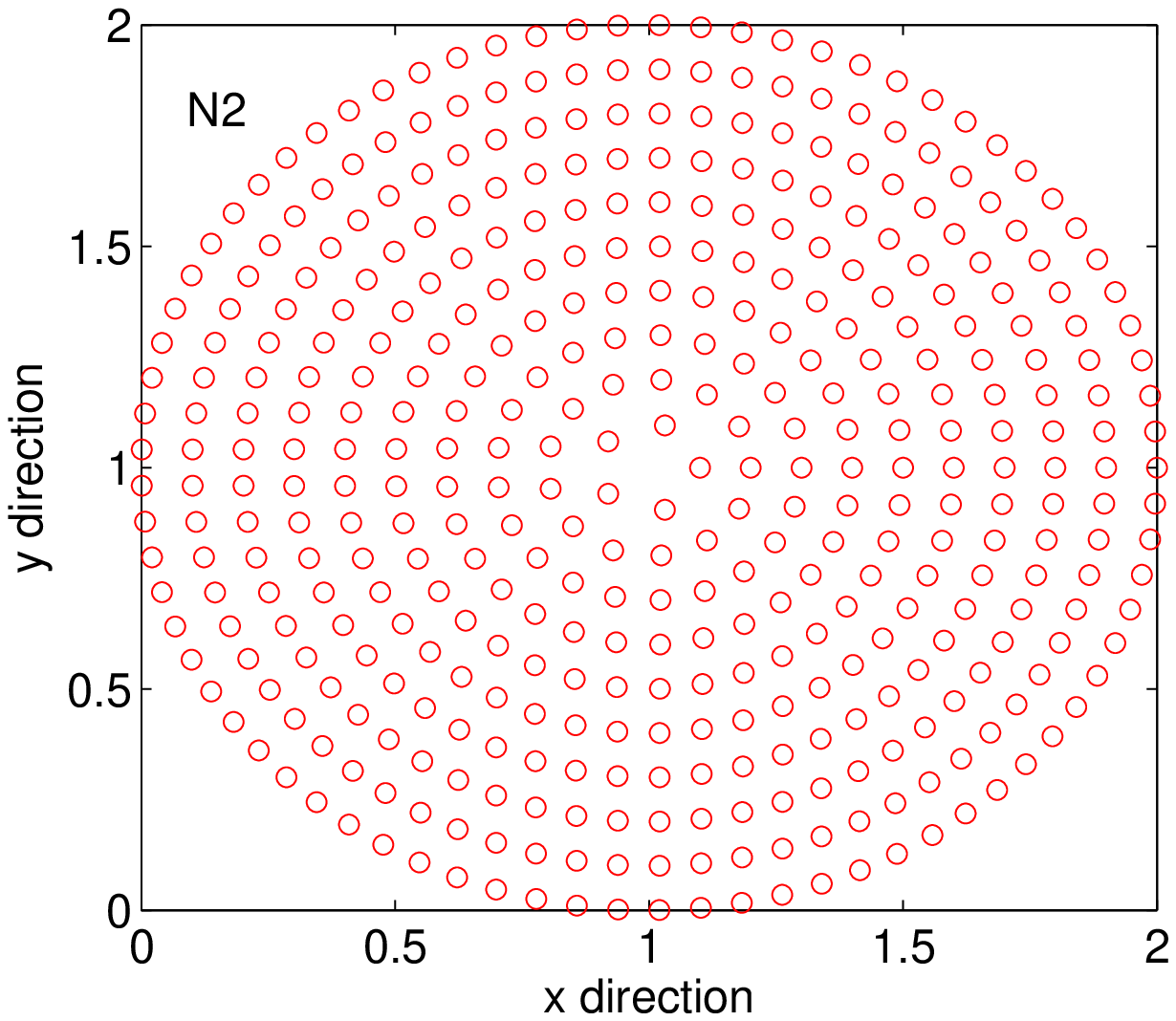}
\label{fig:forceBalanceComparison_b}}
\subfigure{
\includegraphics[width=0.45\linewidth]{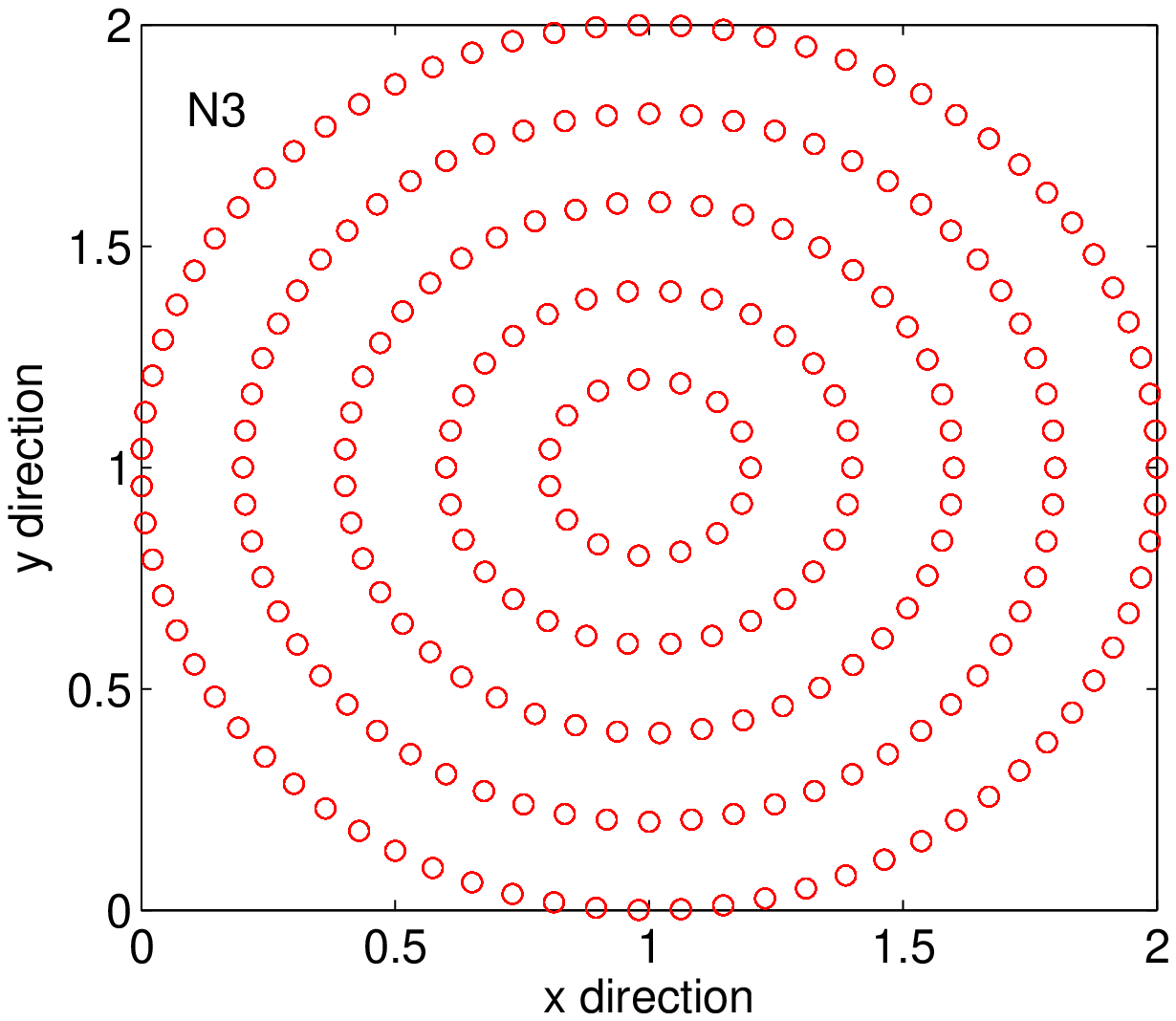}
\label{fig:forceBalanceComparison_b}}
\caption{Three kinds of node collocation approaches. The first one is random distribution (N1), the second one is uniformly distribution in $10$ circles with increment of radius is $0.1$ (N2), the last one is uniformly distribution in $5$ circles with increment of radius is $0.2$ (N3).} \label{fig7}
\end{figure}

Here we also consider the two-dimensional diffusion equation (\ref{eq31}) in a unit disk domain $(x,y)\in \Omega=\sqrt{(x-1)^2+(y-1)^2}\leq 1$. The analytical solution can also be obtained by changing boundary and initial conditions accordingly. To investigate the influence of node collocation methods on the accuracy of numerical results, we adopt three kinds of node collocation approaches presented in Fig. \ref{fig7}. In the first node collocation approach (N1), we adopt a random node-distribution mode. Node collocation approaches N2 and N3 are regular node-distribution modes with different distribution rules. In the approach N2, the collocation nodes are uniformly distributed within $10$ circles having an increment radius of $0.1$; while the collocation nodes in N3 are uniformly distributed in $5$ circles with a radius increment of $0.2$. The main purpose of this is to investigate possible differences between the random and regular node-distribution modes, and determine the influence of node numbers on numerical accuracy.

\begin{table}
\centering
\caption{ABEs ($l_\infty$) of three kinds of collocation approaches for the two-dimensional spatiotemporal FDE with $\alpha=0.7$, $\beta_x=1.6$, and $\beta_y=1.8$ in a unit disk domain at time $t=10$.}
\begin{tabular}{llllll}
\hline\noalign{\smallskip}
N1 & n=140 & n=200 & n=400 & n=500 & c   \\
$l_\infty$ & 0.00350982 & 0.00258652 & 0.00230509 & 0.00239093 & 0.10   \\
\noalign{\smallskip}\hline\noalign{\smallskip}
N2 & n=140 & n=200 & n=400 & n=500 & c   \\
$l_\infty$  & 0.00261043 & 0.00228390 & 0.00201757 & 0.00201826 & 0.15 \\
\noalign{\smallskip}\hline\noalign{\smallskip}
N3 & n=200 & n=300 & n=400 & n=500 & c   \\
$l_\infty$  & 0.00366864 & 0.00366776 & 0.00366775 & 0.00366775 & 0.15 \\
\noalign{\smallskip}\hline
\end{tabular}
\label{tab4}
\end{table}

MAEs of these three collocation approaches are presented in Tab. \ref{tab4}. Generally speaking, the regular node-distribution modes (N2 and N3) are superior to the random node-distribution mode (N1), but the error rates of all three modes are on the same order. It means that the random node-distribution mode (N1), which can be easily obtained for real-world engineering problems, offers acceptable numerical results when compared with the regular node-distribution modes. Meanwhile, the numerical result with a smaller radius increment (e.g., N2 with $\Delta r=0.1$) performed better than the larger radius increment (N3 with $\Delta r=0.2$). However, it should be pointed out that the MAEs of the three approaches become stable when the collocation node number increased. Hence we remark here that a better node-distribution mode is more attractive than simply increasing the collocation node number in the RBF method.

\section{Applications}
\subsection{Two-dimensional spatiotemporal FDE in a circular domain (continuous case)}
\begin{figure}[htb]
\centering
\subfigure{
\includegraphics[width=0.4\linewidth]{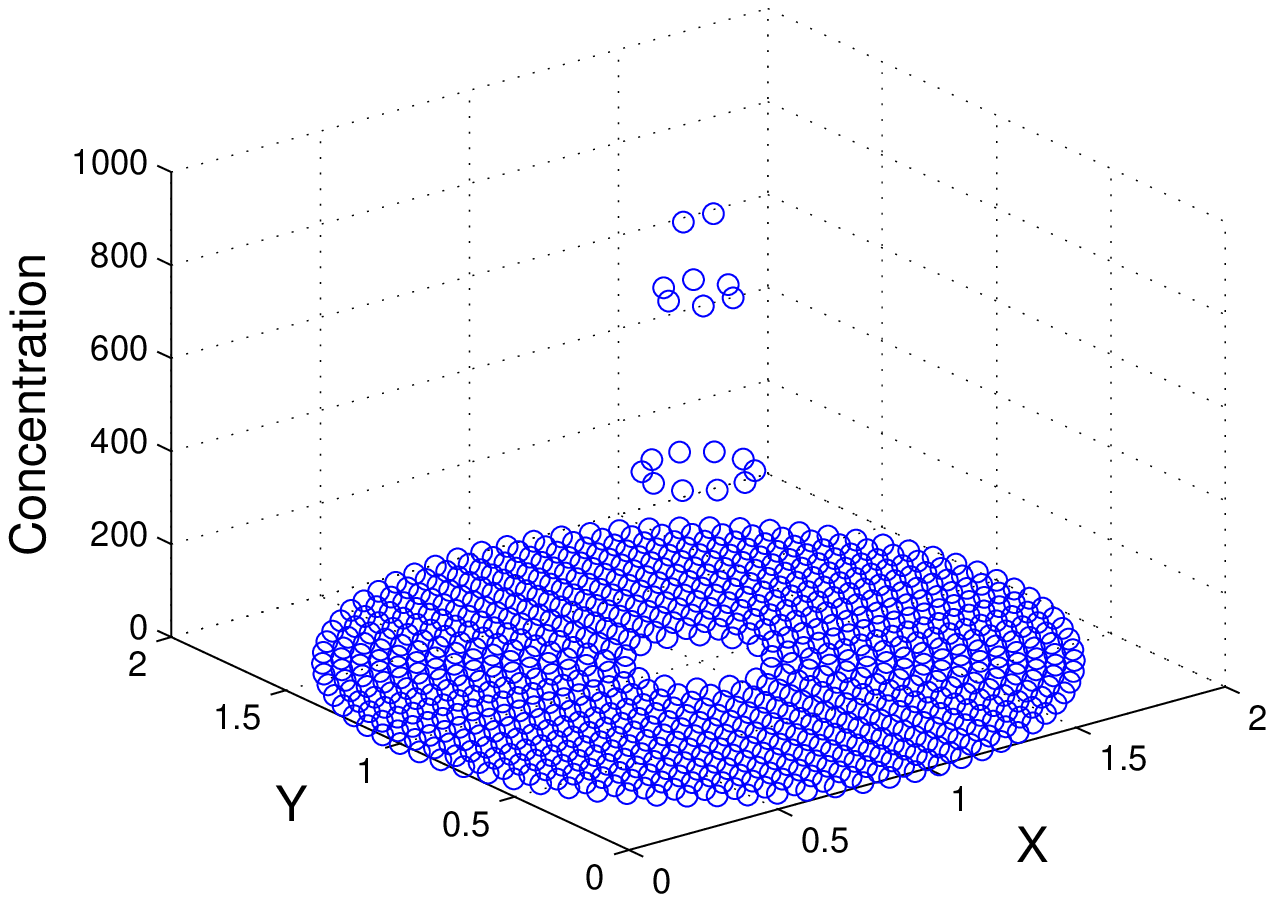}
\label{fig:forceBalanceComparison_a}}
\subfigure{
\includegraphics[width=0.5\linewidth]{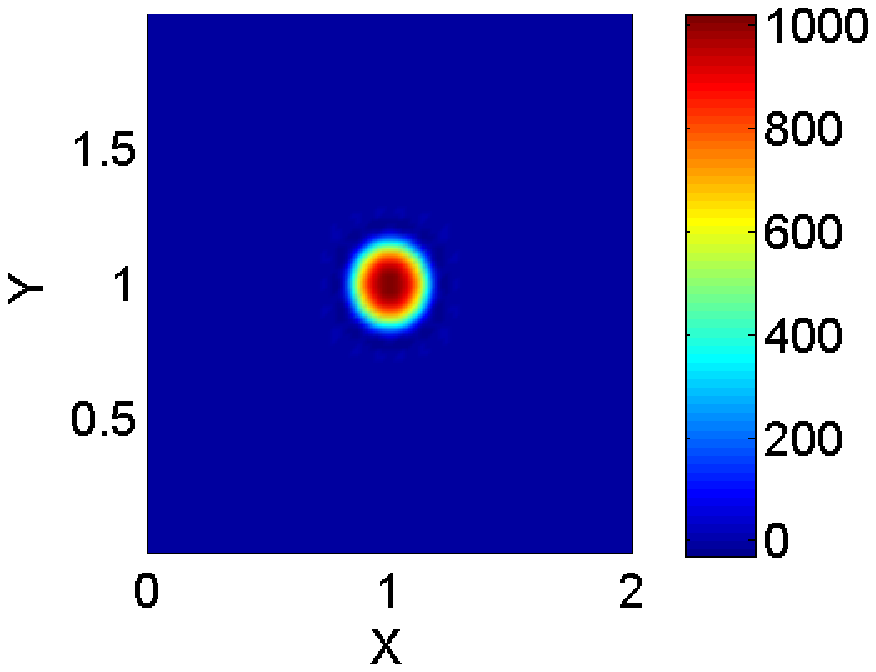}
\label{fig:forceBalanceComparison_b}}
\subfigure{
\includegraphics[width=0.4\linewidth]{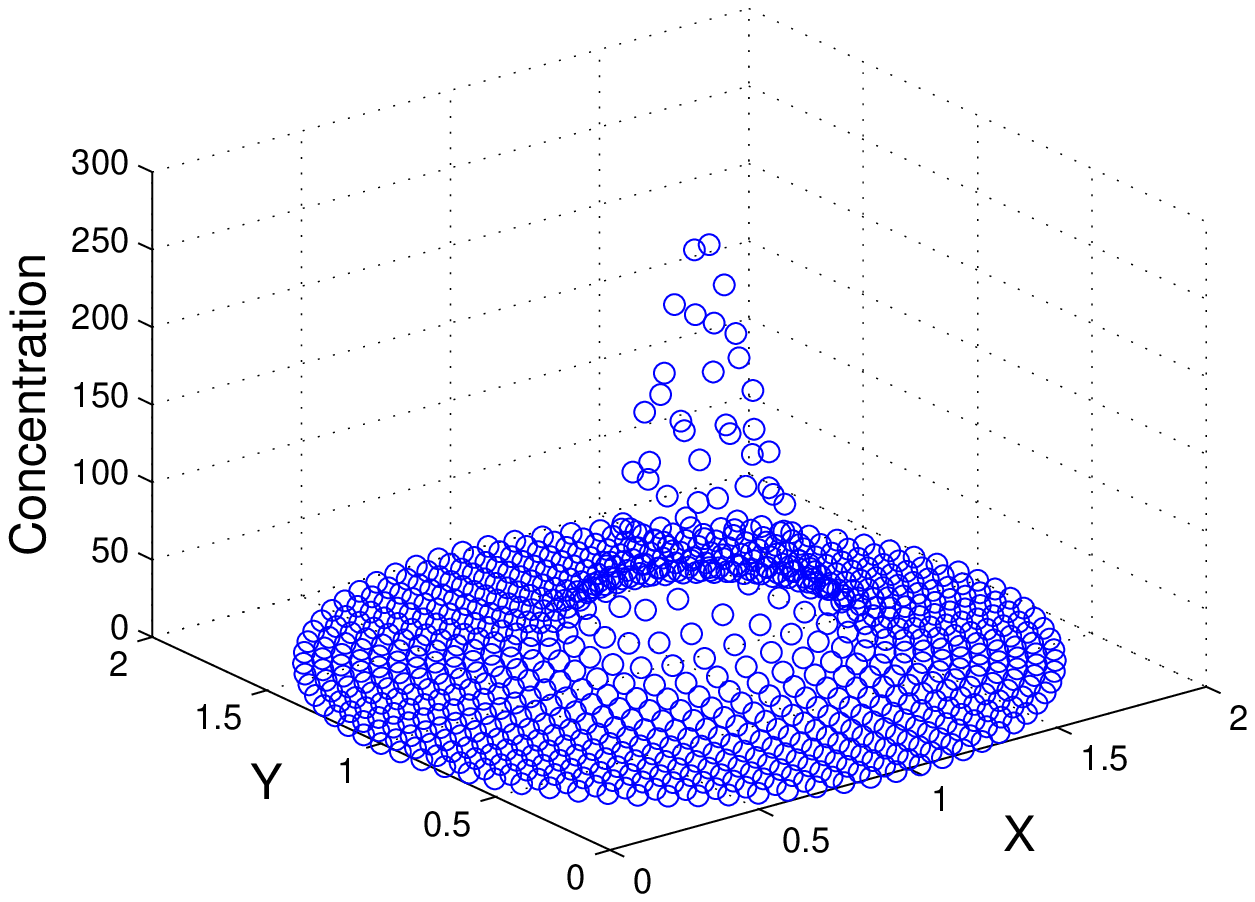}
\label{fig:forceBalanceComparison_a}}
\subfigure{
\includegraphics[width=0.5\linewidth]{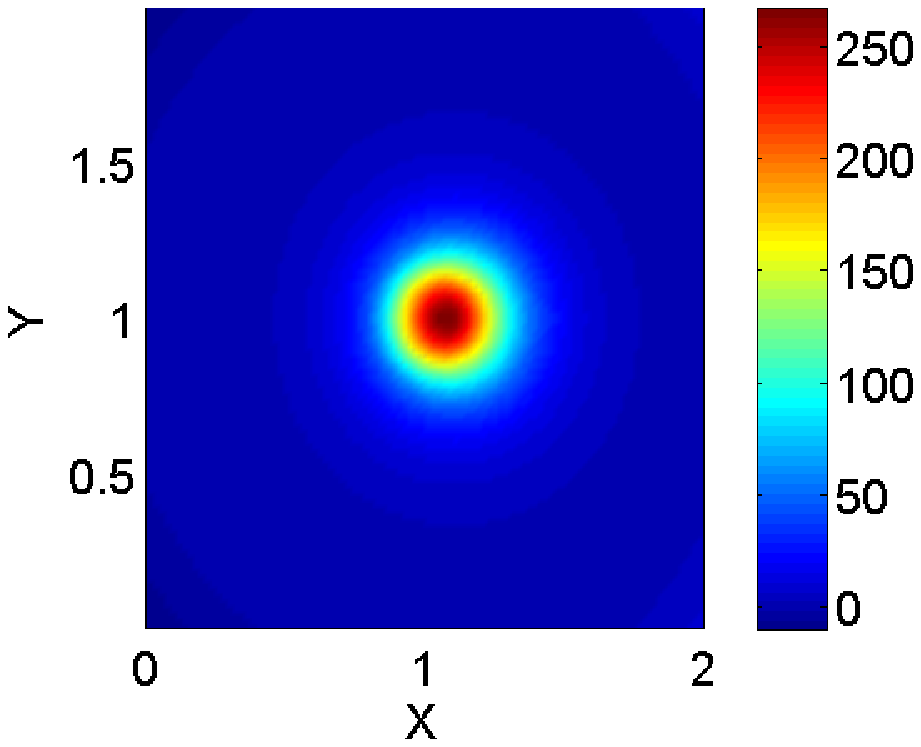}
\label{fig:forceBalanceComparison_b}}
\subfigure{
\includegraphics[width=0.4\linewidth]{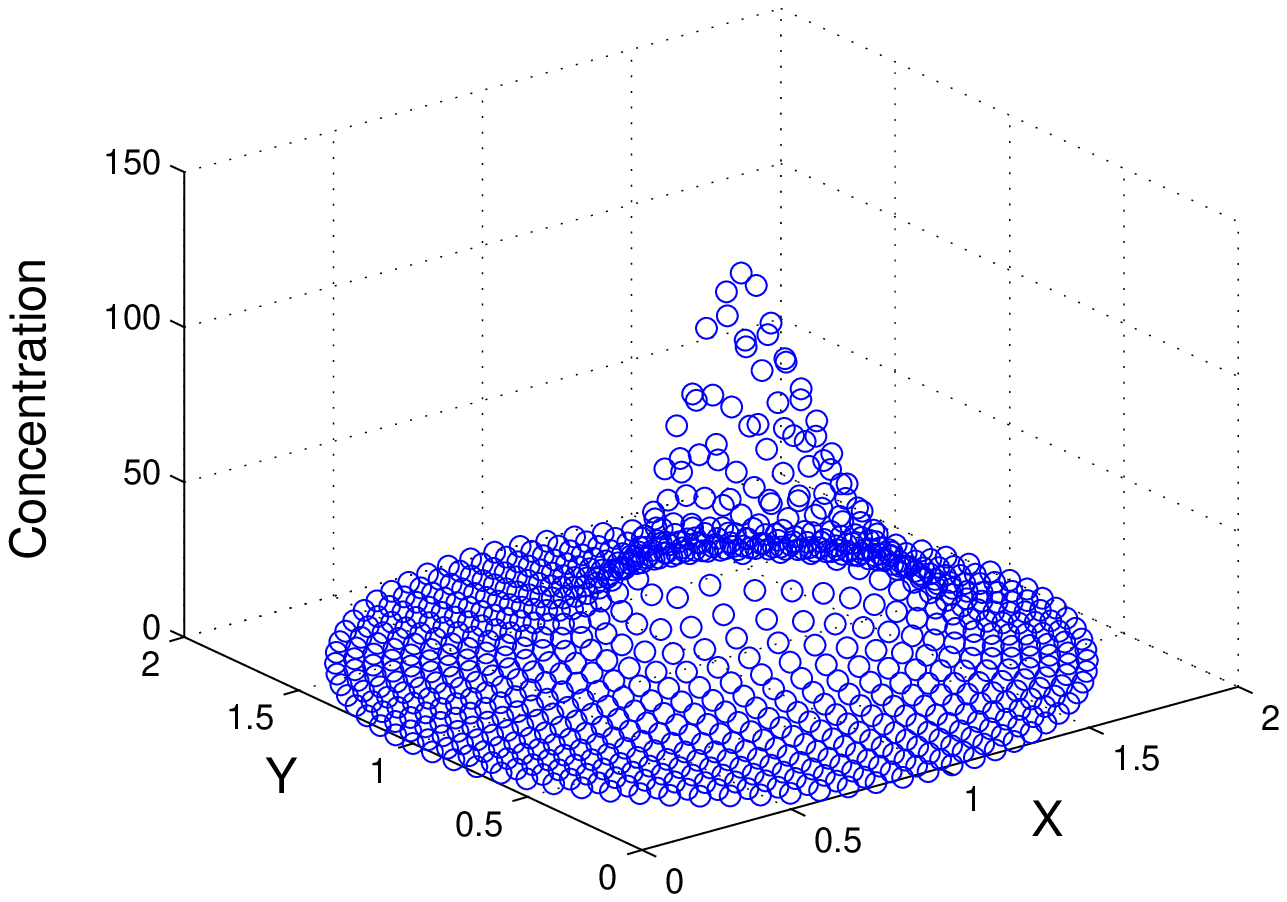}
\label{fig:forceBalanceComparison_a}}
\subfigure{
\includegraphics[width=0.5\linewidth]{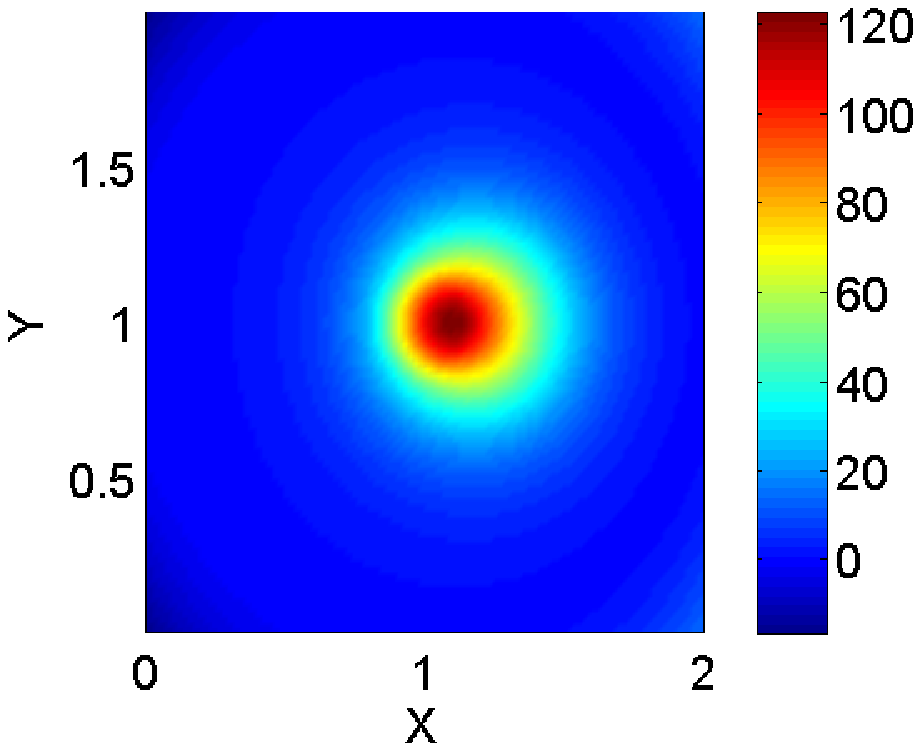}
\label{fig:forceBalanceComparison_b}}
\caption{Concentrations calculated using the proposed scheme for Eq. \ref{eq32} with flow velocity $V=0.012$, dispersion coefficient $K=0.03/2 \pi$, the order of the time fractional derivative $\alpha=0.9$,
and the order of the space fractional derivative $\beta=1.1$. In the numerical scheme, the node number $N=800$, and the shape parameter $C=0.2$. The concentration evolution contour maps for $t=0,\,\,10,\,\,20$ are drawn respectively from top to bottom, with the pointed diagram (left) and the colored nephogram(right). In the colored nephogram, the concentration outside the circular domain remains zero.} \label{fig8}
\end{figure}

\begin{figure}
 \includegraphics[width=0.8\textwidth]{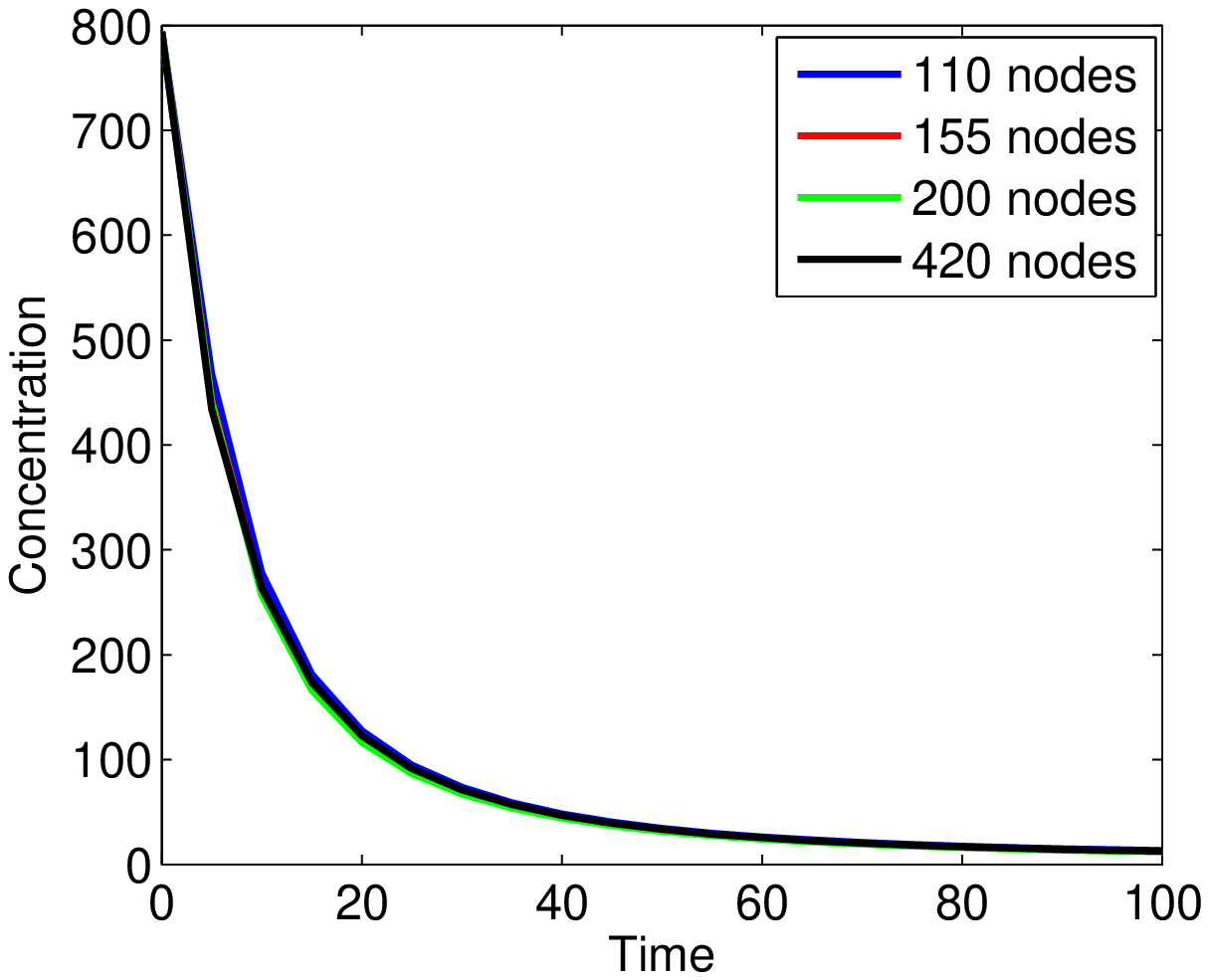}
\caption{Concentration evolution curves computed by different node numbers at (x=1.1, y=1.0). The dimensionless parameters in model Eq. (\ref{eq32}) are: flow velocity $V=0.012$, dispersion coefficient $K=0.03/2 \pi$, the order of the time fractional derivative $\alpha=0.9$, and the order of the space fractional derivative $\beta=1.1$. In the numerical scheme, the shape parameter $C=0.2$. }
\label{fig9}       % Give a unique label
\end{figure}

We consider a two-dimensional spatiotemporal FDE with a compactly supported initial-value condition
\begin{eqnarray}
\left\{
\begin{array}{c}
\displaystyle{\frac{\partial^\alpha u}{\partial t^\alpha} =-V  \frac{\partial u}{\partial x}+K \int_0^{2 \pi} D_\theta^\beta u d \theta,}\\
\displaystyle{(x,y)\in \Omega=(x-1)^2+(y-1)^2\leq 1,\,\,0<\alpha \leq 1,\,\, 1<\beta\leq 2,}\\
\displaystyle{u|_{\partial \Omega} =0,}\\
u(x,y,0)=\left\{
\begin{array}{c}
\displaystyle{1000*2^{1-(1-(x-1)^2/0.04-(y-1)^2/0.04)^{-1}},}\\
\displaystyle{\mbox{if}\,\, (x-1)^2+(y-1)^2<0.04,}\\
\displaystyle{0,\,\, \mbox{otherwise}.}
\end{array}
\right.\;
\end{array}
\right.\;
\label{eq32}
\end{eqnarray}
This model can simulate contaminant transport in isotropic geological media.
The continuous-type space fractional derivative describes similar behaviors for anomalous dispersion along all directions. From the physical point of view, a global correlation of contaminant transport can be better described using the above model than with the previous fractional derivative model which considers non-locality only along $x$ and $y$ coordinates. The initial condition of the above model means that contaminant particles are injected into an isotropic medium at a circular
region $(x,y)\in \Omega=(x-1)^2+(y-1)^2<0.04$. To preserve a homogeneous boundary condition,
here we let $u|_{\partial \Omega} =0$, which means contaminants are not allowed to reach the boundary
of the computational domain during the observation period. We employ the node collocation approach N2 with 800 nodes to calculate the concentration field in the circular domain. The corresponding contour plots of the concentration field for Eq. (\ref{eq32}) are shown in Fig.\ref{fig8}. Fig. \ref{fig8}
illustrates that the proposed scheme is capable of describing the spreading and drifting behaviors of a contaminant plume governed by the model (\ref{eq32}). Since the analytical solution of Eq. (\ref{eq32}) cannot be achieved directly, we compare the numerical concentration evolution curves computed with different node numbers at a fixed point ($x=1.1, y=1.0$). Clearly, a stable result of the numerical concentration evolution curves with different node numbers are obtained from the observation of Fig.\ref{fig9}.

\subsection{Two-dimensional spatiotemporal FDE (discrete case)}%Ö÷³ÌÐòexample5
\begin{figure}[htb]
\centering
\subfigure{
\includegraphics[width=0.45\linewidth]{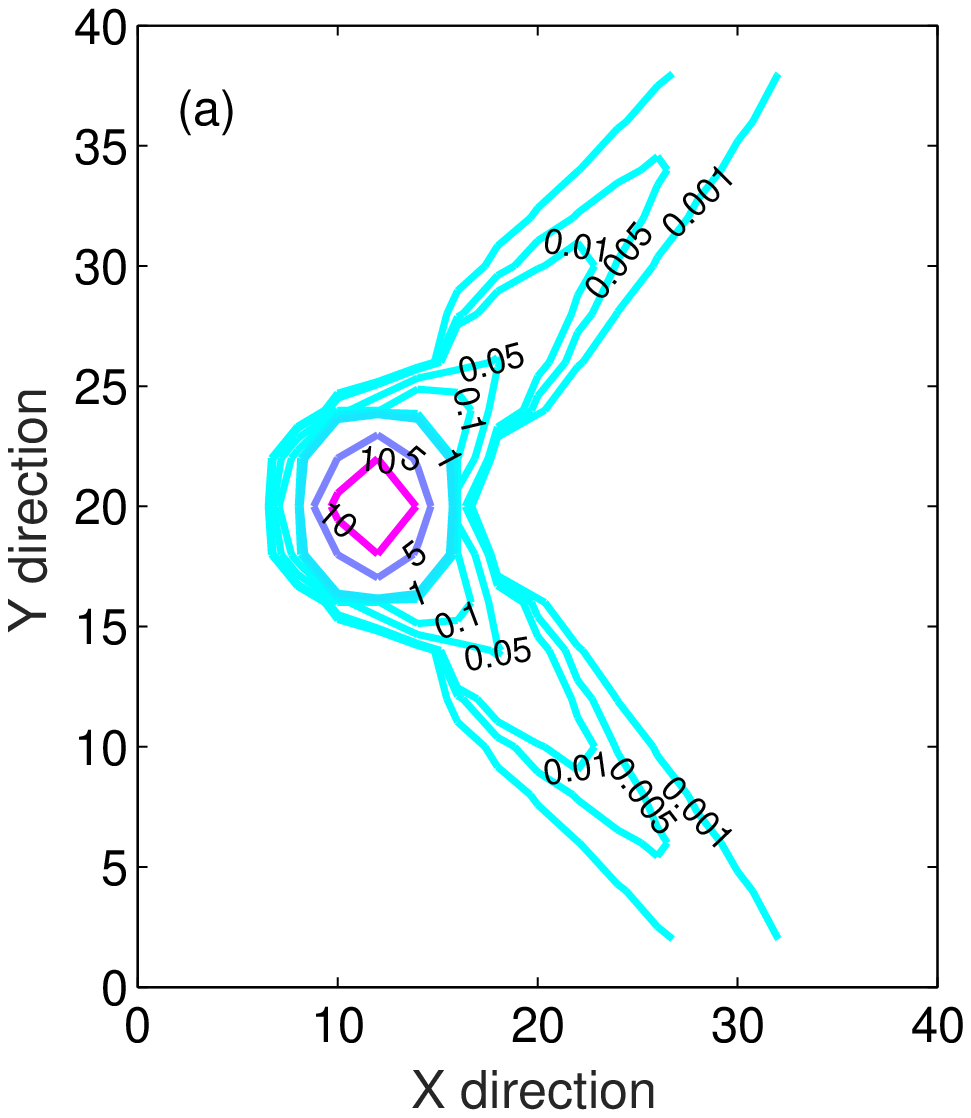}
\label{fig:forceBalanceComparison_a}}
\subfigure{
\includegraphics[width=0.45\linewidth]{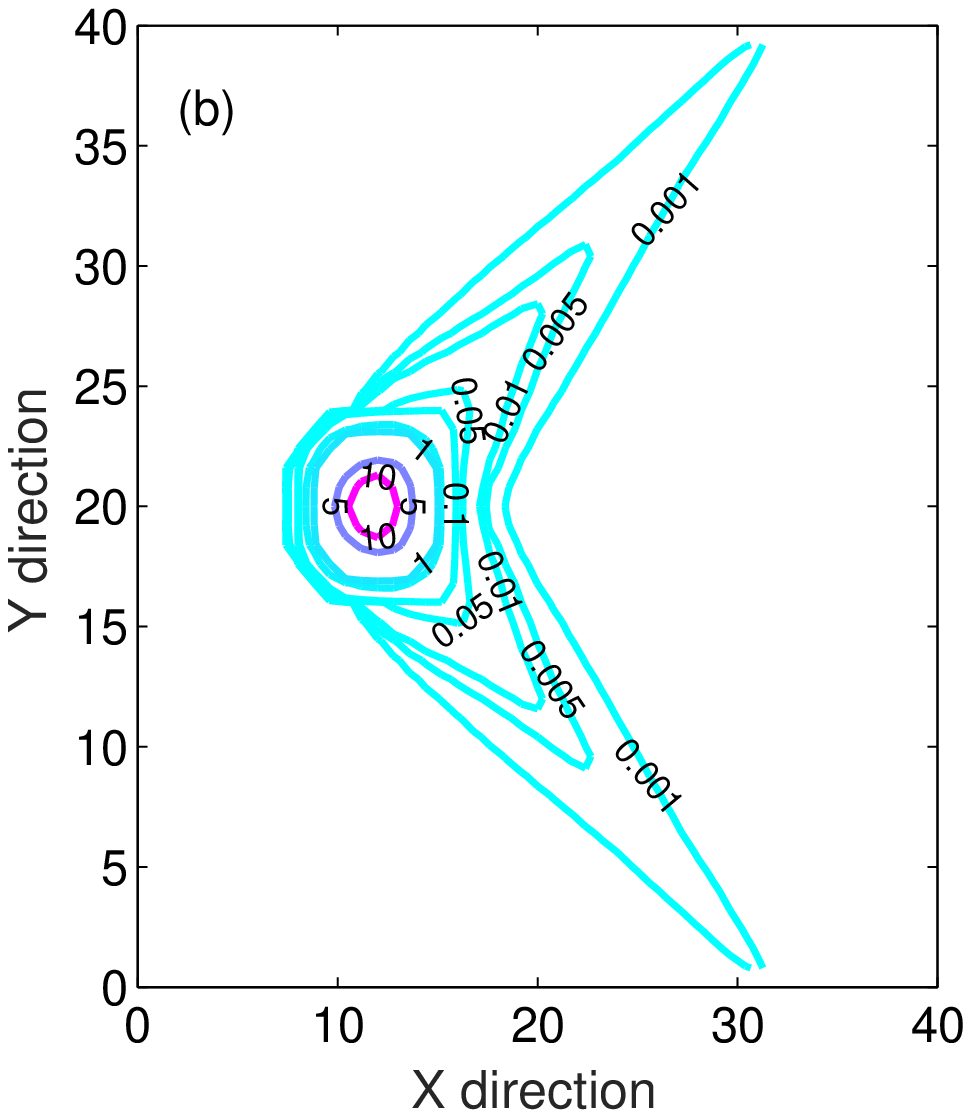}
\label{fig:forceBalanceComparison_b}}
\subfigure{
\includegraphics[width=0.45\linewidth]{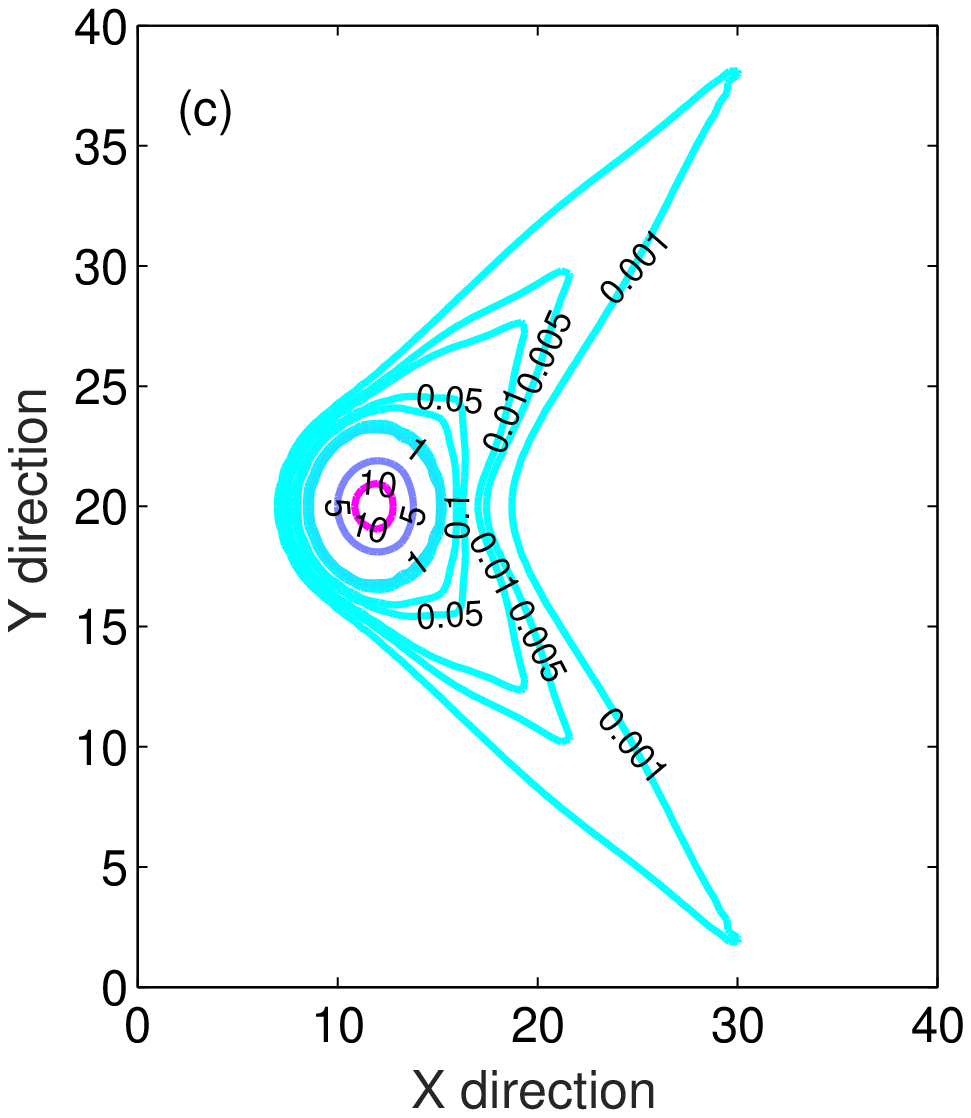}
\label{fig:forceBalanceComparison_a}}
\subfigure{
\includegraphics[width=0.45\linewidth]{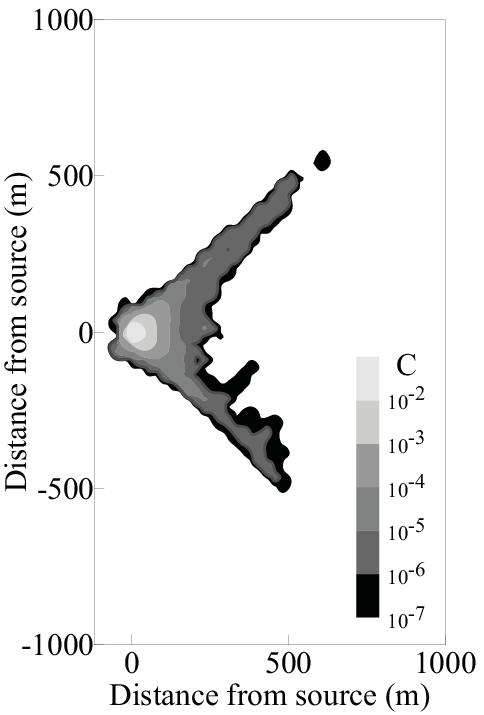}
\label{fig:forceBalanceComparison_a}}
\caption{Concentration field calculated by the proposed scheme for Eq. (\ref{eq33}) with different node numbers. From top to bottom, the node number is $N=441$ (a), $2601$ (b) and $10201$ (c).  The shape parameter $C=0.3$. In the model (\ref{eq33}), diffusion coefficients $K_1=K_2=0.1$, the order of the time fractional derivative $\alpha=0.50$, and the order of the space fractional derivative $\beta=1.55$. The last figure is the experimental data obtained by Reeves et al. \cite{Reeves2008}.} \label{fig10}
\end{figure}

In real-world applications, the occurrence of preferential flow path in aquifers makes the medium an anisotropic one, resulting in direction-dependent anomalous dispersion \cite{Schumer2003}. Therefore, numerical simulations of direction-dependent contaminant transport are of practical importance. For example, Schumer et al. \cite{Schumer2003} used a fast Fourier Transform method to compute multiscaling FDEs. Reeves et al. \cite{Reeves2008} employed the Monte-Carlo method to simulate conservative solute transport through random $2$-D, regional-scale, synthetic fracture networks. Here we establish the following discrete type, spatiotemporal FDE as an example to illustrate how to characterize anomalous transport in an anisotropic medium:
\begin{eqnarray}
\left\{
\begin{array}{c}
\displaystyle{\frac{\partial^\alpha u}{\partial t^\alpha} =K_1 \frac{\partial^\beta u}{\partial (r_{\theta=\pi/4})^\beta}+K_2 \frac{\partial^\beta u}{\partial (r_{\theta=7 \pi/4})^\beta},}\\
\displaystyle{r=\sqrt{(x-12)^2+(y-20)^2}\,\, \mbox{and}\,\,(x,y)\in \Omega=(0,40)*(0,40),}\\
u|_{\partial \Omega} =0,\\
u(x,y,0)=\left\{
\begin{array}{c}
\displaystyle{10/(r+0.1),}\,\,\mbox{if}\,\, \displaystyle{3-r>0,}\\
\displaystyle{0,\,\, otherwise,}
\end{array}
\right.\;
\end{array}
\right.\;
\label{eq33}
\end{eqnarray}
where $K_1$ and $K_2$ are the diffusion coefficients along directions $\theta=1 \pi/4$ and $\theta=7 \pi/4$, respectively. With the help of the generalized fractional derivative used in this model, the diffusion directions do not need to be orthogonal, and the model can be applied to irregular noncontinuum fracture networks and aquifers. To test the influence the number of nodes has on numerical accuracy, here we compare the differences resulting from three numbers of nodes in Fig.\ref{fig10}. The comparison results illustrate that more details of solute diffusion behavior, especially at the tail part, can be captured with a greater number of nodes. Here we also show the Monte-Carlo result of plume concentration in fracture networks \cite{Reeves2008}. A comparison implies that the present model may reasonably characterize anomalous diffusive behavior in anisotropic aquifers, as also concluded in previous works by Reeves et al. (Figure 5) \cite{Reeves2008} and Schumer et al. (Figure 10) \cite{Schumer2003}.

%Èý·ùͼÔÚͼÖÐÈ¡µã·Ö±ðΪ21*21,51*51,101*101

\section{Discussions}

Although the above tests consider only the Dirichlet boundary condition, the numerical scheme proposed by this study may be extended to problems with Neumann or mixed boundary conditions. For example, if we consider a Neumann boundary condition in the $x$ direction, then the expression of the boundary condition in Eq. (\ref{eq11}) can be rewritten as:
\begin{equation}\label{eq34}
\displaystyle{0=-\frac{\partial \mathbf{\Phi}_b}{\partial x}\lambda+\mathbf{u}_b,}
\end{equation}
and the term $\Phi_b$ in Eq. (\ref{eq16}) should be replaced by $\partial \Phi_b / \partial x$. Finally, we can obtain the same expression as Eq. (\ref{eq18}). Furthermore, we can modify the above numerical scheme by combining the Dirichlet and Neumann boundary conditions, to define a mixed boundary condition.

Another interesting point is that relative errors for long-time computation are slightly smaller than those for short-time estimations, which can be observed in Tab. \ref{tab2} and Fig. \ref{fig9}. This numerical property leads to better predictions of late-time solute transport for the presented method in comparison with other approaches. The main reason may be due to the series expression of the Mittag-Leffler function, in which the errors of historical points are balanced when a large number of historical points are used in a long-time computation.

%different parameter optimizations in contrast to those obtained by the FDM [11]. It is observed from Fig.1 that
%Frank¡¯s and Fasshauer¡¯s parameters will lead the Kansa method to higher accuracy and convergence than the FDM.
%From the failure of Hardy¡¯s optimization, one can see that a good choice of

Shape parameter optimization plays an essential role in the accuracy and convergence rate of the Kansa method. Existing literature provides several parameter optimization methods discussed in our previous study \cite{Pang2015}. In this study, we found that the determination method of the shape parameter $C$ in the Kansa method for spatiotemporal FDEs differs slightly from previous approaches.
More theoretical and numerical analyses are needed to explore this mechanism in a future study.

Generally speaking, the computational cost of numerical methods for spatiotemporal FDEs is affected by two aspects: building a coefficient matrix (i.e., difference, stiffness or interpolation matrix) and solving the resulting linear system. Since we employ an analytical approach to solve the resulting linear system with a time fractional derivative, the computational cost of this aspect can be minimized. Our previous study also indicated that the computational cost of building a coefficient matrix with the Kansa method is similar to the cost of the FDM and smaller than the cost of the FEM \cite{Sun2013,Chen2010}. The Kansa method is more flexible for irregular-domains and can be easily extended to three dimensional cases due to the dimension-independent nature of the RBF. However, the computational cost of the dense discretization matrix remains a critical problem in the presented method. A future investigation in fast matrix-vector multiplication techniques, such as the adaptive cross approximation \cite{Bebendorf2003}, the fast multipole method \cite{Coulier2016}, and other approaches \cite{Wang2012,Zhao2016}, is still needed for real-world applications involving large spatial domains.

Since most previous literature employs the commonly used definition of the space fractional derivative (considering only the fractional derivative in $x$ and $y$ coordinates) in numerical computation of two-dimensional FDEs (e.g., Example 2 in a rectangular domain), here we compare numerical results with the commonly used, vector fractional derivative definitions. Clearly, the two definitions are the same when we only consider vector-type fractional derivatives in the directions of $\theta=0$ and $\theta=\pi/2$ from a theoretical approach, and hence we only explore difference in the numerical computation aspect. The maximum relative errors presented in Tab. \ref{tab5} show that they have the same error rate and that the numerical result for the commonly used definition is only slightly improved over the vector version. This implies that the presented scheme which is based on a vector-space fractional derivative, may be used as a fast alternative to the commonly used definition when solving FDE mdoels.

\begin{table}
\centering
\caption{Maximum relative errors (MRE) of the presented scheme for Example 2 in a rectangular domain at time $t=10$ using the vector (C1) and the classical (C2) space fractional derivatives. In Eq. (\ref{eq31}), the dispersion coefficients $K_x=K_y=-1\Gamma (1.5-\beta)(x+1)^{\beta_x}/(2\Gamma (1.5))$, the order of the time fractional derivative $\alpha=0.7$, and the order of the space fractional derivatives $\beta_x=1.6$ and $\beta_x=1.8$. In the numerical scheme, the MQ parameter $c=0.01$.}
\begin{tabular}{lll}
\hline\noalign{\smallskip}
$\Delta=\Delta x=\Delta y$ & MRE (C1) & MRE (C2)\\
\noalign{\smallskip}\hline\noalign{\smallskip}
$\Delta=1/10$ & 0.01262140 & 0.01611325   \\
$\Delta=1/15$ & 0.00759885 & 0.00723610   \\
$\Delta=1/20$ & 0.00744080 & 0.00469782 \\
$\Delta=1/25$ & 0.00653728 & 0.00439643 \\
\noalign{\smallskip}\hline
\end{tabular}
\label{tab5}
\end{table}
%These issues require further mathematical and numerical analysis as well as algorithmic development in the future.

\section{Conclusions}

This study develops a semi-discrete Kansa method to approximate two-dimensional spatiotemporal FDE models for anomalous transport, where the long-time range prediction in large irregular domains is desirable for many environmental protection applications. In the present scheme, the Kansa method is fist used to discretize a large (irregular) spatial domain without the burden of tedious mesh generation, and then an analytical approach is used to obtain a fast solution for the resultant linear equation system with the time fractional derivative.

Numerical examples have shown that the semi-discrete Kansa method can obtain accurate estimations with relatively low computational cost, even using a random distribution node mode. The numerical results of the two-dimensional spatiotemporal FDE with either a continuous or discrete weight function show that the solver has application potential in simulating non-Fickian diffusion for contaminants in real-world, irregular and anisotropic media. Further applications of the proposed method will be contributed in a systematic analysis of the experimental data and a comprehensive comparison with Monte-Carlo simulations \cite{Boggs1993}.

{\bf Acknowledgment.}  We thank Prof. Wen Chen for valuable discussions on Kansa method. The work was supported by the National Natural Science Foundation of China (Grant Nos. 11572112, 41628202, and 11528205). Y. Zhang was also partially supported by the National Science Foundation grant DMS-1460319
and the University of Alabama. This paper does not necessarily reflect the view of the funding agency. %We also thank anonymous reviewers for suggestions that
%significantly improved the presentation of this work.

%\bibliographystyle{plain}

\end{document}